\documentclass[11pt]{article}
\usepackage{latexsym}
\usepackage{fullpage}
\usepackage{amssymb}
\usepackage{amscd}
\def\be{\begin{equation}}
\def\ee{\end{equation}}
\def\bea{\begin{eqnarray*}}
\def\eea{\end{eqnarray*}}
\def\<#1,#2>{\langle\,#1,\,#2\,\rangle}

\def\ra{\rangle}
\def\la{\langle}
\def\R{\mathbb{R}}
\def\C{\mathbb{C}}
\def\H{\mathbb{H}}
\def\RP#1{\ifmmode\R P^{#1}\else$\R P^{#1}$\fi}
\def\CP#1{\ifmmode\C P^{#1}\else$\C P^{#1}$\fi}
\def\HP#1{\ifmmode\H P^{#1}\else$\H P^{#1}$\fi}
\def\SO{\mathrm{SO}}
\def\Sp{\mathrm{Sp}}
\def\Spin{\mathrm{Spin}}
\def\SU{\mathrm{SU}}
\def\U{\mathrm{U}}
\def\G{\mathrm{G}}

\def\id{{\rm id}\,}

\def\pr{{\rm pr}}
\def\tr{{\rm tr}}
\def\Ric{{\rm Ric\,}}

\def\<#1,#2>{\langle#1,#2\rangle}

\newlength{\WidthOfWedge}
\newlength{\WidthOfCirc}
\newlength{\owedgeOffset}

\def\1{\mathbf{1}}
\def\2{\mathit{1}}

\def\Lie{{\mathcal L\,}}

\def\#{\sharp}

\def\es{\,\lrcorner\;}
\def\we{\,\wedge\,}

\def\tfrac#1#2{{\textstyle\frac{#1}{#2}}}
\def\qed{\ifmmode\quad\Box\else$\quad\Box$\fi}
\def\proof{{\bf Proof.}\hskip2mm}
\newtheorem{Pro}{Proposition}[section]
\newtheorem{Lem}[Pro]{Lemma}
\newtheorem{The}[Pro]{Theorem}
\newtheorem{Cor}[Pro]{Corollary}
\newtheorem{Definition}[Pro]{Definition}
\def\beq{\begin{eqnarray*}}
\def\eeq{\end{eqnarray*}}

\newcommand{\firstline}[2]{\makebox[#1][l]{$\displaystyle{#2}$}&&}
\begin{document}
\title{Conformal Killing forms on Riemannian manifolds}
\author{Uwe Semmelmann}
\maketitle
\begin{abstract}
Conformal Killing forms are a natural generalization of conformal vector fields
on Riemannian manifolds. They are defined as sections in the kernel of a
conformally invariant first order differential operator. We show the existence
of conformal Killing forms on nearly K\"ahler and weak $G_2$-manifolds. Moreover,
we give a complete description of special conformal Killing forms. A further
result is a sharp upper bound on the dimension of the space of conformal Killing
forms.
\end{abstract}
%

\section{Introduction}

A classical object of differential geometry are {\it Killing vector fields}.
These are by definition infinitesimal isometries, i.e. the flow of
such a vector field preserves a given metric. The space of all 
Killing vector fields forms the Lie algebra of the isometry
group of a Riemannian manifold and the number of linearly independent
Killing vector fields measures the degree of symmetry of the manifold.
It is known that this number is bounded from above by the
dimension of the isometry group of the standard sphere and,
on compact manifolds, equality is attained if and only if the manifold is isometric
to the standard sphere or the real projective space. Slightly more generally one can consider
{\it conformal vector fields}, i.e. vector fields with a flow
preserving a given conformal class of metrics. There are several geometric
conditions which force a conformal vector field to be Killing.

Much less is known about a rather natural
generalization of conformal vector fields, the so-called 
{\it conformal Killing forms}. These
are p-forms $\psi$ satisfying for any vector field
$X$ the differential equation
\begin{equation}\label{erstedef}
\nabla_X\,\psi
\;-\;
\tfrac{1}{  p+1}\,X\es d\psi 
\;+\;
\tfrac{1}{  n-p+1}\, X^*\,\wedge\,d^*\psi
\;=\;0\ ,
\end{equation}
where $n$ is the dimension
of the manifold, $\nabla$ denotes the covariant
derivative of the Levi-Civita connection, $X^*$ is 1-form dual
to $X$ and  $\es$ is the operation dual
to the wedge product. It is easy to see that a conformal Killing
1-form is dual to a conformal vector field. Coclosed conformal
Killing $p$-forms are called {\it Killing forms}. For $p=1$ 
they are dual to Killing vector fields. 

The left hand side
of equation~(\ref{erstedef}) defines a first order elliptic differential
operator $T$, which was already studied in the context of Stein-Weiss
operators (c.f.~\cite{branson}). Equivalently one can describe a conformal 
Killing form as
a form in the kernel of $T$. From this point of view
conformal Killing forms are  similar to twistor spinors in spin
geometry. One shared property is the conformal invariance of
the defining equation. In particular, any form which is parallel
for some metric $g$, and thus a Killing form for trivial reasons,
induces non-parallel conformal Killing forms for  metrics
conformally equivalent to $g$ (by a non-trivial change of the metric). 

Killing forms, as a generalization of the Killing vector fields, 
were introduced by K.~Yano in \cite{yano}. Later
S.~Tachibana (c.f.~\cite{ta3}), for the case of 2--forms,
and more generally T.~Kashiwada (c.f.~\cite{ka}, \cite{ta4}) introduced
conformal Killing forms generalizing conformal vector fields.

Already K.~Yano noted that a $p$--form $\psi$ is a Killing form
if and only if for any geodesic $\gamma$ the $(p-1)$--form 
$\,{\dot \gamma} \,\lrcorner\, \psi\,$ is parallel along $\gamma$. 
In particular,
Killing forms give rise to quadratic first integrals of the
geodesic equation, i.e. functions which are constant along geodesics. 
Hence, they can be used to integrate the equation
of motion. This was first done in the article~\cite{penrose}
of R.Penrose and M. Walker, which  initiated an intense study of
Killing forms in the physics literature. In particular, there is
a local classification of Lorentz manifolds with Killing 2-forms.
More recently Killing forms and conformal Killing forms have been
successfully applied to define symmetries of field equations
(c.f.~\cite{be1}, \cite{be2}).

Despite this longstanding interest in (conformal) Killing forms 
there are only very few global results on Riemannian manifolds. 
Moreover the number of the known non-trivial examples on compact 
manifolds is surprisingly small. The aim of this article is to
fill this gap and to start a study of global properties of 
conformal Killing forms. 

As a first contribution we will show that there are several classes of
Riemannian manifolds admitting Killing forms, which so far did not appear
in the literature. In particular, we will show
that there are Killing forms on nearly K{\"a}hler manifolds and
on manifolds with a weak $G_2$--structure. All these examples are
related to Killing spinors and nearly parallel vector cross products.
Moreover, they are all so-called {\it special Killing forms}. The restriction
from Killing forms to special Killing forms is analogous to the
definition of a Sasakian structure as a unit length Killing vector
field satisfying an additional equation. One of our main results
in this paper is the complete description of manifolds
admitting special Killing forms.

Since conformal Killing forms are sections in the kernel of an
elliptic operator  it is clear that they 
span a finite dimensional space in the case of compact manifolds.
Our second main result is an explicit upper bound for the dimension 
of the space of conformal Killing forms on arbitrary connected
Riemannian manifolds. The upper bound is provided by the dimension
of the corresponding space on the standard sphere. It is also shown 
that if the upper bound is attained the manifold has to be conformally flat.

In our paper we  tried to collect all that is presently known for 
conformal Killing forms on Riemannian manifolds. This includes some
new proofs and new versions of known results.

\bigskip

{\bf Acknowledgments}

In the first place, I would like to thank Prof.~D.~Kotschick
for  valuable discussions, his support and his interest in my work.
I am grateful to A.~Moroianu and G.~Weingart for
many helpful comments,  important hints and a continued interest in
the topic of conformal Killing forms on Riemannian manifolds.


\section{The definition of conformal Killing forms}\label{definition}

In this section, we will introduce conformal Killing forms, give integrability
conditions and state several well-known elementary properties,
including equivalent characterizations.

Let $( V,\,\la\cdot,\cdot\ra)$ be an $n$--dimensional Euclidean vector space. 
Then the $O(n)$--representa\-tion  $V^*\otimes\Lambda^pV^*$ has the following 
decomposition:
\begin{equation}\label{deco1}
V^*\otimes\Lambda^pV^*
\;\cong\;
\Lambda^{p-1}V^* \oplus \Lambda^{p+1}V^* \oplus \Lambda^{p,1}V^* \ ,
\end{equation}
where $\Lambda^{p,1}V^*$ is the intersection of the kernels of wedge
 product and 
contraction map. The highest weight of the representation 
$\Lambda^{p,1}V^*$ is the sum of the highest weights of $V^*$ and of 
$\Lambda^pV^*$.
Elements of $\,\Lambda^{p,1}V^* \subset V^* \otimes \Lambda^pV^*\,$ 
can be
considered as 1-forms on $V$ with values in $\,\Lambda^pV^*$.
For any $v\in V$, $\alpha \in V^*$ and $\psi \in \Lambda^pV^*$,
the projection 
$\pr_{\Lambda^{p,1}}:V^*\otimes\Lambda^pV^*\, \rightarrow 
\Lambda^{p,1}V^*\,$ is then explicitly given by
\be\label{prtwist}
[\pr_{\Lambda^{p,1}} (\alpha\, \otimes\, \psi)]\, v
\;\;:=\;\;
\alpha (v)\, \psi 
\;-\;
\tfrac{1}{ p+1}\, v\es(\alpha\,\wedge\,\psi)
\;-\;
\tfrac{1}{ n-p+1}\, v^* \we (\alpha^\#\es\psi) \ ,
\ee
where $v^*$ denotes the 1-form dual to $v$, i.e. $v^*(w) = \la v, w\ra$,
$\alpha^\#$ is the vector defined by $\alpha(v) = \la \alpha^\#, v\ra$
and $ v \es$ denotes the interior multiplication which is dual to the
wedge product $v \we$.

This decomposition immediately translates to Riemannian manifolds
$\,(M^n,\,g)$, where we have the decomposition
\be\label{deco}
T^*M\otimes\Lambda^pT^*M
\;\cong\;
\Lambda^{p-1}T^*M
\oplus
\Lambda^{p+1}T^*M
\oplus
\Lambda^{p,1}T^*M
\ee
with $\Lambda^{p,1}T^*M$ denoting the vector bundle corresponding to 
the representation $\Lambda^{p,1}$. The covariant derivative $\nabla \psi$ 
of a $p$--form $\psi$ 
is a section of $\;T^*M\otimes\Lambda^pT^*M$, projecting it onto
the summands $\,\Lambda^{p+1}T^*M\,$ resp. $\,\Lambda^{p-1}T^*M\,$ yields $d\psi$ resp. $d^*\psi$.
The projection onto the third summand $\,\Lambda^{p,1}T^*M\,$ defines a
natural first order differential operator $T$, which we will call the {\it twistor operator}. The twistor operator
$ 
T:\Gamma(\Lambda^p T^*M)\,\rightarrow \,\Gamma(\Lambda^{p,1}T^*M) 
\subset 
\Gamma(T^*M\otimes\Lambda^pT^*M)
$
is given for any vector field $X$ by the following formula
$$
[\,T\psi\,]\,(X)
\;:=\;
[\pr_{\Lambda^{p,1}}(\nabla \psi)]\,( X)
\;=\;
\nabla_X\, \psi
\;-\;
\tfrac{1}{  p+1}\, X \es d\psi
\;+\;
\tfrac{1}{ n-p+1}\, X^*\,\wedge\,d^*\psi \ .
$$
This definition  is similar to the definition of the twistor operator 
in spin geometry, where one has the decomposition of the
tensor product of spinor bundle and cotangent bundle into the sum
of spinor bundle and kernel of the Clifford multiplication. The 
twistor operator is defined as the projection of the covariant derivative
of a spinor onto the kernel of the Clifford multiplication, which,
as a vector bundle, is  associated to the representation
given by the sum of highest weights of spin and standard representation.


\begin{Definition}
A p-form $\,\psi$ is called a {\it conformal Killing p-form} if and only if 
$\,\psi\,$ is in the kernel of $\,T$,~i.e. 
if and only if $\,\psi\,$ satisfies for all vector fields $X$ the equation
\begin{equation}\label{killing}
\nabla_X\,\psi\;=\;
\tfrac{1}{  p+1}\,X\es d\psi \;-\;
\tfrac{1}{  n-p+1}\, X^*\,\wedge\,d^*\psi\ .
\end{equation}
If the p-form $\psi$ is in addition coclosed it is called a 
{\it Killing p-form}. 
This is equivalent to $\nabla\psi\in\Gamma(\Lambda^{p+1}T^*M)$ or to 
$X \es \nabla_X \psi = 0$ for any vector field $X$. 
\end{Definition}

Closed conformal
Killing forms  will be called $\ast$--{\it Killing forms}.
In the physics literature, equation~(\ref{killing}) defining a conformal
Killing form is often called the {\it Killing--Yano equation}. A
further natural notation is {\it twistor forms}, which is also motivated
by the following observation.
Let $(M, g)$ be a spin manifold and let $\psi$ be a {\it twistor spinor},~i.e. 
a section of the spinor bundle lying in the kernel of the spinorial twistor
operator or, equivalently, a spinor satisfying for all vector fields $X$ the 
equation
$
\nabla_X \psi \;=\; -\,\frac{1}{n}\,X \cdot D\psi
$,
where $D$ is the Dirac operator and $\,\cdot\,$ denotes the Clifford multiplication. Given two such twistor spinors,
$\psi_1$ and $\psi_2$, we can introduce k-forms $\omega_k$, which are
on any tangent vectors $X_1, \ldots , X_k$ defined by
$$
\omega_k(X_1, \ldots, X_k)
:=
\la (X_1^* \wedge \ldots \wedge X_k^*)\cdot \psi_1, \psi_2 \ra \ .
$$
It is well-known that for $k=1$ the form $\omega_1$ is dual to
a conformal vector field. Moreover, if $\psi_1$ and $\psi_2$ are
Killing spinors the form $\omega_1$ is dual to a Killing
vector field. Recall that  a Killing spinor is a section $\psi$ of the
spinor bundle satisfying for all vector fields $X$ and some constant $c$
the equation $\nabla_X \psi = c X \cdot \psi$,~i.e. Killing spinors are
special solutions of the twistor equation. More generally we have
\begin{Pro}\label{twistor}
Let $(M^n, g)$ be a Riemannian spin manifold with twistor 
spinors $\,\psi_1$ and $\psi_2$. Then for any $k$ the associated
k-form $\omega_k$ is a conformal Killing form. 
\end{Pro}

The proof, which follows from a simple local calculation, is given in the
appendix. 

Decomposition (\ref{deco}) implies that the covariant derivative 
$\nabla \psi$ splits into three components. Using the twistor operator $T$ 
we can write the covariant derivative of a $p$-form $\psi$ as
\be\label{covariant}
\nabla_X\psi\;=\;
\tfrac{1}{  p+1}\,X\es d\psi \;-\;
\tfrac{1}{  n-p+1}\, X^*\,\wedge\,d^*\psi\;+\;
[\,T\psi\,]\,(X) \ .
\ee
This formula leads to the following pointwise norm estimate together with
a further characterization of conformal Killing forms (c.f.~\cite{galot}).
\begin{Lem}\label{kernel}
Let $(M^n,\,g)$ be a Riemannian manifold and let $\,\psi\,$ be any $p$--form. 
Then
\be\label{kernel1}
| \,\nabla \psi \,|^2
\;\ge\;
\tfrac{1}{  p+1}\,| \,d\psi\,|^2\;+\;\tfrac{1}{  n-p+1}\,|\,d^*\psi\,|^2 \ ,
\ee
with equality if and only if $\psi$ is a conformal Killing $p$--form.
\end{Lem}
As an application of Lemma \ref{kernel} one can prove that the Hodge 
star-operator $\ast$ maps conformal Killing $p$--forms into conformal Killing 
$(n-p)$--forms. In particular, $\ast$ interchanges closed and coclosed 
conformal Killing form.

Differentiating equation~(\ref{covariant}) we obtain two Weitzenb{\"o}ck
formulas, which play an important role in the proof of many global
results. Similar characterizations were obtained in~\cite{ka}. 
For any p-form $\psi$ we have the equations
\begin{eqnarray}
\nabla^*\nabla\psi 
&=& \label{weiz1}
\tfrac{1}{ p+1}\,
d^*d\, \psi\;\;+\;\; \tfrac{1}{  n-p+1}\,d d^*\, \psi\;\;+\;\;T^*T\,
\psi \ ,
\\[1.5ex]
q(R)\,  \psi 
&=&\label{weiz2}
\tfrac{p}{ p+1}\,d^*d\, \psi\;\;+\;\;\tfrac{n-p}{ n-p+1}\,d d^*\, \psi
\;\;-\;\;T^*T\, \psi \ ,
\end{eqnarray}
where $q(R)$ is the curvature term appearing in the classical
Weitzenb{\"o}ck formula for the Laplacian on $p$-forms:
$\Delta  = d^*d\,+\,dd^* = \nabla^*\nabla \,+\, q(R) $. It is  the symmetric
endomorphism of the bundle of differential forms defined by
$$
q(R)\;=\;\sum\,e^*_j\,\wedge\,e_i \es R_{e_i,e_j},
$$
where $\,\{e_i\}\,$ is any local ortho-normal frame and $R_{e_i,e_j}$
denotes the curvature of the form bundle.
On forms of degree one and two one has an explicit expression for the action of
$q(R)$. Indeed, if $\xi$ is any 1--form, then $q(R)\,\xi = \Ric(\xi)$ and
if $\omega$ is any 2--form 
then
\be\label{curvop}
q(R)\,\omega
\;=\;
\Ric(\omega) -2{\cal R} \,(\omega),
\ee
where $\Ric$ denotes the symmetric endomorphism of the form bundle
obtained by extending the Ricci curvature as derivation. Moreover,
$\;{\cal R}$ denotes the Riemannian curvature operator defined 
on vector fields $\, X, Y, Z, U \, $ by
$\,  g({\cal R}(X \we Y),\, Z \we U) = - g(R(X,\,Y)\,Z,\,U)$.

Integrating the second Weitzenb{\"o}ck formula (\ref{weiz2}) gives rise 
to an important integrability condition and a characterization
of conformal Killing forms on compact manifolds. Indeed we have
\begin{Pro}\label{integrabl}
Let $(M^n,\,g)$ a compact Riemannian manifold. Then a p--form is
a conformal Killing $p$--form, if and only if
\be\label{char}
q(R)\,  \psi 
\;\;=\;\;
\tfrac{p}{ p+1}\,d^*d\, \psi\;\;+\;\;\tfrac{n-p}{ n-p+1}\,d d^*\, \psi \ .
\ee
\end{Pro}
This proposition implies that there are no
conformal Killing forms on compact manifolds, where $q(R)$ has
only negative eigenvalues. This is the case on manifolds with constant
negative sectional curvature or on conformally flat manifolds with
negative-definite Ricci tensor. For coclosed forms, Proposition~\ref{integrabl} 
is a  generalization of the well-known characterization of Killing vector 
fields on compact manifolds, as divergence free vector fields in the 
kernel of $\,\Delta - 2\,\Ric$. In the general case, it can be 
reformulated as
\begin{Cor}\label{integrab}
Let $(M^n,\,g)$ a compact Riemannian manifold with a coclosed $p$--form
$\psi$. Then $\psi$ is a Killing form if and only if
$$
\Delta \, \psi
\;\;=\;\;
{p+1\over p}\,q(R)\,  \psi \ .
$$
\end{Cor}
We note that there are similar results for $\ast$--Killing forms 
and for conformal Killing m-forms on 2m-dimensional manifolds.

A further interesting property of the equation defining
conformal Killing forms is its conformal invariance (c.f.~\cite{be2}).
The precise formulation  is
\begin{Pro}\label{conformal}
Let $(M^n,\,g)$ be a Riemannian manifold with a conformal Killing
p--form $\psi$. Then ${\widehat \psi}:= e^{(p+1)\lambda}\psi$ is a
conformal Killing p--form with respect to the conformally equivalent
metric ${\hat g}:= e^{2\lambda}g$.
\end{Pro}
In particular, it follows from this proposition that the Lie derivative
with respect to conformal vector fields preserves the space of
conformal Killing forms.

There is still another characterization of conformal Killing 
forms which is often given as the definition. 

\begin{Pro}\label{killdef}
Let $(M^n,\,g)$ be a Riemannian manifold. A $p$--form $\psi$
is a conformal Killing form if and only if there exists a  $(p-1)$--form
$\,\theta\,$ such that
\bea
\firstline{.5cm}{(\nabla_Y \,\psi)(X, X_2,\ldots,X_p)
\;\; + \;\;
(\nabla_X \,\psi)(Y, X_2,\ldots,X_p)}\\[1.5ex]
&& = 
2 g(X,\,Y)\,\theta (X_2,\ldots,X_p)
\,-\,
\sum^p_{a=2}\,(-1)^a
\left(
g(Y,\,X_a)\,\theta (X,\,X_2,\ldots,{\hat X}_a,\ldots,\,X_p)
\right.\\[1.5ex]
&&
\phantom{2\,g(X,\,Y)\,\theta }
\left.
\;\;+\;\;g(X,\,X_a)\,\theta (Y,\,X_2,\ldots,{\hat X}_a,\ldots,\,X_p)
\right) 
\eea
for any vector fields $Y, X, X_1,\ldots X_p$, where ${\hat X}_a$
means that $X_a$ is omitted.
\end{Pro}

It was already mentioned in the introduction that the interest
in Killing forms in relativity theory stems from the fact that they
define first integrals of the geodesic equation. At the end of this 
chapter, we will now describe this construction in more detail.
Let $\,\psi\,$ be a Killing $p$-form and let $\,\gamma\,$ be a 
geodesic,~i.e.
$\,\nabla_{\,\dot \gamma} \, {\dot \gamma} = 0$. Then
$$
\nabla_{\dot \gamma} \,({\dot \gamma} \es \psi)
\;=\;
(\nabla_{\dot \gamma} \, {\dot \gamma}) \es \psi
\;+\;
{\dot \gamma} \es \nabla_{\dot \gamma} \, \psi
\;=\; 0 \ ,
$$
i.e. ${\,\dot \gamma} \es \psi\,$ is a $(p-1)$--form parallel along
the geodesic $\gamma$ and in particular its length is constant
along $\gamma$. The definition of this constant can be given in a
more general context. Indeed for any p-form $\psi$ we can
consider a symmetric bilinear form $K_\psi$ defined 
for any vector fields $X, Y$ as
$$
K_\psi(X,\,Y) \;:=\; g( X \es \psi,\, Y \es \psi\,)\ .
$$

\noindent
For Killing forms  the associated bilinear form has a very nice property.
\begin{Lem}
If $\psi$ is a Killing form, then the associated symmetric bilinear form 
$K_\psi$ is a Killing tensor,~i.e. for any vector fields $X, Y, Z$
it satisfies 
the equation
\begin{equation}\label{killing2}
(\nabla_X K_\psi)(Y,\,Z) \;+\; (\nabla_Y K_\psi)(Z,\,X) 
  \;+\; (\nabla_Z K_\psi)(X,\,Y) \;=\; 0 \ .
\end{equation}
In particular, $K_\psi ({\dot \gamma},\,{\dot \gamma})$
is constant along any geodesic $\gamma$.
\end{Lem}

In general, a $\,(0,k)$--tensor $\,\mathcal T\,$ is called 
{\it Killing tensor} if the
complete symmetrization of $\,\nabla {\mathcal T}\,$ vanishes. 
This is equivalent
to $\,(\nabla_X {\mathcal T})(X,\ldots,X) = 0$. It follows again that 
for such a
Killing tensor, the expression $\,{\mathcal T}({\dot \gamma},\ldots,
{\dot \gamma})\,$ is 
constant along any geodesic $\,\gamma\,$ and hence defines a $\,k$-th 
order
first integral of the geodesic equation. Note that the length of the
$\,(p-1)$--form $\,X \es \psi\,$ is $\,K_\psi(X, X)\,$ and that 
$\,\tr(K_\psi) = p \,|\psi|^2$.


\section{Examples of conformal Killing forms}\label{examples}
\setcounter{equation}{0}

We start with parallel forms which are obviously in the kernel of the 
twistor operator and thus are conformal Killing forms. 
Using Proposition \ref{conformal}, we see that with any parallel 
p-form $\,\psi$, the form 
$\,{\widehat \psi}:= e^{(p+1)\lambda}\,\psi\,$ is a conformal Killing 
$\,p$--form with respect to the conformally equivalent
metric $\,{\widehat g}:= e^{2\lambda}\,g$. This new form $\,{\widehat \psi}\,$
is in general no longer parallel. 

\medskip

Conformal Killing forms were introduced as a generalization
of conformal vector fields, i.e. we have the following 
well-known result.

\begin{Pro}
Let $(M,\,g)$ be a Riemannian manifold. Then a vector field $\,\xi$ is dual
to a conformal Killing 1-form if and only if it  
is a conformal vector field,~i.e. if there exists a 
function $f$ such that $\Lie_\xi \,g = f \, g$. Moreover, $\xi$ is dual to a Killing 
1-form if and only if it is a Killing vector field,~i.e. 
if $\Lie_\xi \,g = 0$.
\end{Pro}

The simplest examples of manifolds with conformal Killing forms
are the spaces of constant curvature. We will recall the result for the 
standard sphere $(S^n,\,g)$ with scalar curvature
$s=n(n-1)$. The spectrum of the Laplace operator on p-forms consists of
two series: 
$$
\lambda_k'
\;=\;(p+k)(n-p+k+1)
\qquad \qquad
\mbox{and}
\qquad \qquad
\lambda_k''
\;=\;(p+k+1)(n-p+k) \ ,
$$
where $k=0,1,2,\ldots$.
The eigenvalues $\lambda_k'$ correspond to closed eigenforms, whereas the
eigenvalues $\lambda_k''\,$ correspond to coclosed eigenforms. The 
multiplicities
of the eigenvalues are well--known. In particular, we find for the minimal 
eigenvalues $\lambda_0'$ and $\lambda_0''\,$ that
$$
\lambda_0' 
\quad \mbox{has multiplicity}\quad
{n+1 \choose p}
\qquad
\mbox{and}
\qquad
\lambda_0''
\quad \mbox{has multiplicity}\quad
{n+1 \choose p+1} \ .
$$

\noindent
The conformal
Killing forms turn out to be sums of eigenforms of the Laplacian 
corresponding to 
the minimal eigenvalues on $\, \ker(d)\,$ resp. $\,\ker(d^*)\,$. 

\begin{Pro}\label{sphere}
A p--form $\omega$ on the standard sphere $(S^n,\,g)$ is a conformal
Killing form, if and only if it is a sum  of eigenforms for  the 
eigenvalue $\,\lambda_0'\,$ resp. 
of eigenforms for the eigenvalue $\,\lambda_0''$.
\end{Pro}

The first interesting class of manifolds admitting conformal
Killing forms are Sasakian manifolds. These are contact
manifolds satisfying a normality (or integrability) condition.
In the context of conformal Killing forms, it is convenient to
use the following 

\noindent
\begin{Definition} 
A Riemannian manifold $(M,\,g)$ is called a {\it Sasakian} manifold,
if there exists a unit length Killing vector field $\,\xi\,$ satisfying
for any vector field $X$ the equation
\begin{equation}\label{sasaki1}
\nabla_X\,(d\,\xi^*) \;=\; -\,2\,X^* \wedge \; \xi^* \ .
\end{equation}
\end{Definition}

Note that in the usual definition of a Sasakian structure, as a
special contact structure one has the additional condition 
$\,\phi^2 = -\,\id\,+\,\eta \otimes\xi\,$ for the associated
endomorphism $\phi=-\nabla\xi$ and the 1-form $\eta := \xi^*$. 
But this equation is implied by~(\ref{sasaki1}), if we write~(\ref{sasaki1})  
first as 
\be\label{sasaki2}
(\nabla_X\phi)(Y)\;=\;g(X,\,Y)\,\xi\;-\;\eta(Y)\,X \ ,
\ee
and take then the scalar product with $\xi$. It  follows
that the dimension of a Sasakian manifolds has to be odd and
if $\,\dim (M) = 2n+1$, then $\,\xi^*\wedge\,(d\,\xi^*)^n$ is
the Riemannian volume form on $M$.

There are many examples of Sasakian manifolds,~e.g. given as 
$S^1$--bundles over K{\"a}hler manifolds. Even in the special
case of {\it 3-Sasakian} manifolds, where one has three
unit length Killing vector fields, defining Sasakian
structures satisfying the $SO(3)$-commutator relations, one knows
that there are infinitely many diffeomorphism types (c.f.~\cite{galicki}).

On a manifold with a Killing vector field $\xi$ we have the
Killing 1-form $\xi^*$. It is then natural to ask whether
$\,d\,\xi^*\,$ is also a conformal Killing form. The next proposition
shows that for Einstein manifolds this is the case, if and only if 
$\,\xi\,$ defines a Sasakian structure. Slightly more general, we have
\begin{Pro}\label{sasaki}
Let $(M,\,g)$ be a Riemannian manifold with a Sasakian structure
defined by a unit length vector field $\,\xi$. Then the
2--form $\,d \, \xi^*$ is a conformal Killing form. Moreover, if
$(M^n,\,g)$ is an Einstein manifold with scalar curvature $\,s$ 
normalized to $s=n(n-1)$ and if $\,\xi\,$ is a unit length Killing
vector field such that $\,d\,\xi^*\,$ is a conformal Killing form,
then $\,\xi\,$ defines a Sasakian structure.
\end{Pro}
\proof
We first prove that for a Killing vector field $\,\xi\,$ defining
a Sasakian structure, the 2-form  $\,d \xi^*\,$ is a conformal Killing form.
From the definition~(\ref{sasaki1}) of the Sasakian structure  we obtain:
$\, d^* d \,\xi^* = 2(n-1)\,\xi^*$. Substituting $\,\xi^*\,$ in (\ref{sasaki1})
using this formula yields:
$$
\nabla_X\,(\,d\,\xi^*) \;=\; -\, 2 X^* \we \tfrac{1}{2(n-1)}\,d^* d \,\xi^*
\;=\; -\tfrac{1}{n-1}\,X^* \we d^* d\, \xi^* \ .
$$
But since $\,d\,\xi^*$, is closed  this equation implies that $\,d\,\xi^*\,$
is indeed a conformal Killing form. 
To prove the second statement, we first note that
$\,d^*d\,\xi^* = \Delta \,\xi^* = 2 \Ric \,(\xi^*) = 2(n-1)\xi^*\,$ because
of equation (\ref{char}) for Killing 1-forms and the assumption that 
$(M,\,g)$ is an Einstein manifold with normalized scalar curvature. 
Then we can reformulate the condition that $\,d \xi^*\,$ is a closed
conformal Killing form to obtain
$$
\nabla_X\,(\,d\,\xi^*) \;=\; -\tfrac{1}{n-1}\,X^* \we d^* d\, \xi^*
\;=\;-\tfrac{1}{n-1}\,X^* \we 2(n-1)\xi^*
\;=\;-2\,X^* \we \xi^* \ ,
$$
i.e. the unit length Killing vector field $\, \xi\,$ also satisfies
the equation (\ref{sasaki1}) and thus defines a Sasakian structure.
\qed

We know already that on a Sasakian manifold defined by a Killing vector
field $\,\xi$, the dual 1-form $\,\xi^*\,$ and the 2-form $\,d\,\xi^*\,$ are
both conformal Killing forms. In fact, the same is true for all possible
wedge products of $\,\xi^*\,$ and $\,d\,\xi^*$. We have
\begin{Pro}\label{powers}
Let $(M^{2n+1},\,g,\,\xi)$ be a Sasakian manifold with Killing vector field
$\,\xi$. Then 
$$
\omega_k 
\;\;:=\;\; 
\xi^* \we (d\xi^*)^k
$$
is a Killing $\,(2k+1)$--form for $\,k=0,\ldots n$. Moreover, $\,\omega_k\,$
satisfies for any vector field $\,X\,$ and any $\,k\,$ the additional equation
$$
\nabla_X (d\omega_k)
\;\;=\;\;
-\,2\,(k+1)\, X^* \we \omega_k .
$$
In particular, $\,\omega_k\,$ is an eigenform of the Laplace operator 
corresponding 
to the eigenvalue $\,4(k+1)(n-k)$.
\end{Pro}
This can be proved by a simple local calculation. However, it is also
part of a more general property which we will further discuss in 
Section~\ref{special}.

Recall that a form $\psi$ on a Sasakian manifold is called {\it horizontal} if
$\xi \es \psi = 0$, where $\xi$ is the vector field defining the
Sasakian structure. In~\cite{yama3.5} resp. \cite{yama4} S.~Yamaguchi
proved the following

\begin{The}
Let $(M,\,g)$ be a compact Sasakian manifold, then
\begin{enumerate}
\item
any horizontal conformal Killing form of odd degree is Killing, and
\item
any conformal Killing form of even degree has a 
unique decomposition \\
into the sum of a Killing form and a 
$\ast$--Killing form.
\end{enumerate}
\end{The}

\medskip

We will now describe a general construction which provides new examples of
Killing forms in degrees 2 and 3. For this aim we have to recall
the notion of a vector cross product (c.f.~\cite{gray}). Let $V$ be a 
finite dimensional
real vector space and let $\la\cdot,\cdot\ra$ be a non-degenerate
bilinear form on $V$. Then a {\it vector cross product} on $V$ is 
defined as a 
linear map $P:V^{\otimes r}\rightarrow V$ satisfying the axioms
\bea
(i)
\quad \la P(v_1,\,\ldots,\,v_r),\,v_i\ra
&=&
0 \qquad (1\le i \le r),\\[1.5ex]
(ii)
\qquad |P(v_1,\,\ldots,\,v_r)|^2
&=&
\det(\la v_i,\,v_j\ra) \ .
\eea
Vector cross products are completely classified. There are only
four possible types: 1-fold and (n-1)-fold vector cross products
on  n-dimensional vector spaces, 2-fold vector cross products on
7-dimensional vector spaces and 3-fold vector cross products on
8-dimensional vector spaces. We will consider r-fold vector cross 
products on Riemannian manifolds $(M,\,g)$. These are tensor fields of
type $(r,\,1)$ which are fibrewise r-fold vector cross products.
As a special class, one has the so-called {\it nearly parallel} vector
cross products. By definition they satisfy the differential equation
$$
(\nabla_{X_1} P)(X_1,\,\ldots,\,X_r)\;\;=\;\;0
$$
for any vector fields $X_1,\ldots,X_r$. Together with an r-fold vector
cross product $P$, one has an associated $(r+1)$-form $\,\omega\,$
defined as
$$
\omega(X_1,\,\ldots,\,X_{r+1})\;\;=\;\;
g( P(X_1,\,\ldots,\,X_r),\, X_{r+1}) \ .
$$
The definition of a nearly parallel vector cross product
is obviously equivalent to the condition $\,X \es \nabla_X \,\omega
=0\,$ for the associated form. Hence, we obtain
\begin{Lem}\label{np}
Let $P$ be a nearly parallel r-fold vector cross product with associated
form $\,\omega$. Then $\,\omega\,$ is a Killing $(r+1)$-form.
\end{Lem}

We will examine the four possible types of vector cross
products to see which examples of manifolds with Killing forms
one can obtain. We start with 1-fold vector cross products,
which are equivalent to almost complex structures 
compatible with the metric. Hence, a Riemannian manifold $(M,\,g)$ 
with a nearly parallel 1-fold vector cross product $J$
is the same as an almost Hermitian manifold, where the almost
complex structure $\,J\,$ satisfies $\,(\nabla_X\,J)\,X\,=\,0\,$ for all 
vector fields $\,X$. Such manifolds are called {\it nearly K{\"a}hler}.
It follows from Lemma \ref{np} that the associated 2-form $\omega$
defined by $\,\omega(X,\,Y)=g(JX,\,Y)\,$ is a Killing 2-form. 
On a K{\"a}hler manifold, $\,\omega\,$ is the K{\"a}hler form  and thus parallel
by definition. But there are also many non-K{\"a}hler, nearly K{\"a}hler 
manifolds,~e.g.
the 3-symmetric spaces which were classified by A.~Gray and J.~Wolf
(c.f.~\cite{graywolf}). 
Due to a result of S.~Salamon (c.f.~\cite{sala}) nearly K{\"a}hler, non-K{\"a}hler 
manifolds are never  Riemannian symmetric spaces.

Next, we consider 2-fold vector cross products. They are defined on
7 dimensional Riemannian manifolds and exist, if and only if the 
structure group 
of the underlying manifold $\,M\,$ can be reduced to the group 
$\,G_2\subset O(7)$,
~i.e. if $M\,$ admits a topological $\,G_2$-structure. Riemannian
manifolds with a nearly parallel 2-fold vector cross product are
called {\it weak} $\,G_2$--{\it manifolds}. There are many examples
of homogeneous and non-homogeneous $\,G_2$--manifolds,~e.g. on any 
7-dimensional 3-Sasakian manifold. Here exists a canonically defined  
(additional) Einstein metric which is weak-$G_2$ (c.f.~\cite{fia}). 

Finally, we have to consider the $\,(n-1)$-fold and 3-fold vector cross 
products. But in these cases, results of A.~Gray show that the 
associated forms have to be parallel (c.f.~\cite{gray}). 
Hence, they yield only trivial examples of conformal Killing forms.

We have seen that nearly K{\"a}hler manifolds are special almost Hermitian
manifolds where the K{\"a}hler form $\,\omega$, defined by 
$\,\omega(X, Y) = g(JX, Y)$, is a Killing 2-form. This leads to the natural
question whether there are other almost Hermitian manifolds, where the 
K{\"a}hler form is a conformal Killing form. The following proposition 
gives an answer to this question. 

\begin{Pro} \label{almherm}
Let $\,(M^{2n},\,g,\,J)\,$ be an almost Hermitian
manifold. Then the K{\"a}hler form $\,\omega$ is a conformal Killing
2-form if and only if the manifold is nearly K{\"a}hler or K{\"a}hler.
\end{Pro}
\proof
Let $\,\Lambda\,$ denote the contraction with the 2-form $\,\omega$, i.e.
$\,\Lambda  = \frac{1}{2} \sum Je_i \es e_i \es$. On an almost
Hermitian 
manifold 
(with K{\"a}hler form
$\,\omega$), one has the following well known formulas:
$$
\Lambda(d\omega) \;=\; J(d^*\omega)
\qquad\quad
\mbox{and}
\qquad\quad
d\omega \;=\; (d\omega)_0 \;+\; \tfrac{1}{n-1}\,(Jd^*\omega) \we \omega \ ,
$$
where $\,(d\omega)_0\,$ denotes the effective or primitive part, i.e. the part
of $\,d\omega\,$ in the kernel of $\,\Lambda$.
We will show that if $\,\omega\,$ is a conformal Killing 2-form, then it
has to be coclosed. The defining equation of a Killing 2-form reads
$$
(\nabla_X\,\omega)(A,\,B)
\;=\;
\tfrac{1}{3}\,
d\omega(X,\,A,\,B)
\;-\;
\tfrac{1}{2n-1}\,
\left(
g(X,\,A)\,d^*\omega(B) \;-\; g(X,\,B)\,d^*\omega(A)
\right) \ .
$$
Because $\,\nabla_X J \circ J + J \circ \nabla_X J = 0\,$ we
see that $\,\nabla_X\omega\,$ is an anti-invariant 2-form. Setting
$X=e_i$ and $A=Je_i$ and summing over an ortho-normal basis $\{e_i\}$
we obtain
\bea
-\,d^*\omega(JB)
&=&
\tfrac{1}{3}\,
\sum \, d\omega(e_i,\,Je_i,\,B) 
\;+\; 
\tfrac{1}{2n-1}\,\sum \, g(e_i,\,B)\,d^*\omega(Je_i)\\[1.5ex]
&=&
\tfrac{2}{3}\,\Lambda(d\omega) \;+\; \tfrac{1}{2n-1}\,d^*\omega(JB)\\[1.5ex]
&=&
(\tfrac{1}{2n-1} \;-\; \tfrac{2}{3})\,d^*\omega(JB) \ .
\eea
From this equation follows immediately $\,d^*\omega=0$, i.e. $\omega$
is already a Killing 2-form. But this is equivalent for $(M,g,J)$ to be nearly
K{\"a}hler, where we consider K{\"a}hler manifolds as a special case of
nearly K{\"a}hler manifolds.
\qed 


\section{Special Killing forms}\label{special}
\setcounter{equation}{0}

In \cite{ta2} S.~Tachibana and W.~Yu introduced the notion of
special Killing forms. This definition seemed to be rather restrictive and indeed
the only discussed examples were spaces of constant curvature.
Nevertheless, it turns out that almost all examples of Killing forms 
described in the preceding section are special and we will now show
that there are only a few further examples.

From the following definition it becomes clear that the restriction 
from Killing forms to special Killing forms is analogous to the restriction 
from Killing vector fields to Sasakian structures. 

\begin{Definition}
A {\it special Killing form} is a Killing form $\,\psi\,$ which  for
some constant $\,c\,$ and any vector field $\,X\,$ satisfies 
the additional equation

\begin{equation}\label{defspecial1}
\nabla_X\,(d\psi) \;=\; c\, X^* \we \psi \ .
\end{equation}
\end{Definition}

There is an equivalent version of equation (\ref{defspecial1}), which
gives a definition closer to the original one. Indeed, a special
Killing form can be defined equivalently as a Killing form satisfying
for some (new) constant $\,c\,$ and for any vector fields 
$\,X,\,Y\,$ the equation
\be\label{defspecial2}
\nabla^2_{X,\,Y}\,\psi
\;\;=\;\;
 c \, \left( g(X,\,Y)\,\psi
\;-\;
X \we Y \es \psi \right) \ .
\ee
From equation (\ref{defspecial1}) it follows immediately that special 
Killing $\,p$--forms are eigenforms of the Laplacian corresponding to the 
eigenvalue $\,-c(n-p)$. Hence, on compact manifolds the constant $\,c\,$ has 
to be negative. 

Our first examples of special Killing forms came from Sasakian manifolds.
Here the defining equation~(\ref{sasaki1}) coincides with 
equation~(\ref{defspecial1}) for the constant $\,c=-2$,~i.e. a Killing
vector field $\,\xi\,$ defining a Sasakian structure is dual to a special
Killing 1-form with constant $\,c=-2$. Moreover, we have seen in
Proposition \ref{powers} that on a Sasakian manifold also the forms
$\,\omega_k := \xi^* \wedge (d\xi^*)^k\,$ are special Killing forms.
All other known examples are given in
\begin{Pro}
The following manifolds admit special Killing forms:
\begin{enumerate}
\item
Sasakian manifolds with defining Killing vector field $\,\xi$. Here all the
Killing forms $\,\omega_k = \xi^* \we (d\xi^*)^k\,$ are special
with constant $\,c=-\,2(k+1)$.
\item
Nearly K{\"a}hler non-K{\"a}hler manifolds in dimension 6. Here the associated 
2--form $\,\omega\,$ is special with constant $\,c=-\,\frac{s}{10}\,$ and the
3-form $\,\ast\, d\, \omega\,$ is special with constant   
$\,c=-\,\frac{2s}{15}$,
where $\,s\,$ denotes the scalar curvature.
\item
Weak $\,G_2$--manifolds of scalar curvature $\,s$. Here the associated 
3--form is a special Killing form of constant $\,c=-\,\frac{2\,s}{21}$.
\item
The standard sphere $\,S^n$ of scalar curvature $\,s=n(n-1)$. Here all Killing 
$\,p$--forms,~i.e. all coclosed  minimal eigenforms of the Laplacian 
 are special with constant  $\,c=-\,(p+1)$.
\end{enumerate}
\end{Pro}
Note that weak $\,G_2$--manifolds
and the nearly K{\"a}hler non-K{\"a}hler manifolds in dimension 6
are Einstein manifolds, hence they have constant scalar curvature.
One can easily see that the associated 2--form on a nearly K{\"a}hler manifold 
of dimension different from 6 is never special. 

We will now give a complete description
of compact Riemannian manifolds admitting special Killing forms. It turns out
that a  $\,p$-form  on $\,M\,$  is a special Killing
form if and only if it induces a $\,(p+1)$--form on the metric
cone $\,{\widehat M}\,$ which is parallel. Since the metric cone is 
either flat or
irreducible, the description of special Killing forms is reduced to a
holonomy problem,~i.e. to the question which holonomies admit parallel
forms. This question can be completely answered and the existence of
parallel forms on the cone can be retranslated into
the existence of special geometric structures on the base manifold. 
The result will be that special Killing forms can exist only on Sasakian
manifolds, nearly K{\"a}hler manifolds or weak $\,G_2$--manifolds.
Our approach here is similar to the one of Ch.~B{\"a}r in \cite{baer} which
lead to the classification of Killing spinors.

The metric cone $\,\widehat M\,$ over a Riemannian
manifold $(M,\,g)$ is defined as a warped product,~i.e. 
$\,\widehat M = M\times \R_+\,$
with metric $\,{\widehat g} := r^2g \,+ \,dr^2$. An easy calculation shows
that the Levi-Civita connection on 1--forms is given by
\bea
{\widehat \nabla}_X \,Y^*
&=&
\nabla_X Y^* \;-\; \tfrac{1}{r}\,g(X,\,Y)\,dr
,\qquad
{\widehat \nabla}_X \,dr
\;=\;
r\,X^* \ ,
\\[1.5ex]
{\widehat \nabla}_{\partial_r}\,X^*
&=&
-\,\tfrac{1}{r}\,X^*,
\qquad
{\widehat \nabla}_{\partial_r} \,dr
\;=\;0 \ ,
\eea
where $\,X, Y\,$ are vector fields tangent to $\,M\,$ with $\,g$--dual
1--forms $\,X^*, Y^*$, and where $\,\partial_r\,$ is the radial vector field
on $\widehat M$ with $dr(\partial_r)=1$. From this we immediately
obtain the following useful formulas
$$
{\widehat \nabla}_X \,\psi
\;=\;
\nabla_X\,\psi
\;-\;
\tfrac{1}{r}\,dr \we ( X \es \psi),\qquad
{\widehat \nabla}_{\partial_r} \,\psi
\;=\;
-\,\tfrac{p}{r}\,\psi \ ,
$$
where $\,\psi\,$ is a $\, p$-form on $\,M\,$ considered as $\,p$--form
on $\,\widehat M$. For any $\,p$--form $\,\psi\,$ on $\,M$, we define an
associated $\,(p+1)$--form $\,\widehat \psi\,$ on $\,\widehat M\,$ by
\begin{equation}\label{coneform}
{\widehat \psi}
\;\;:=\;\;
r^p\, dr \we \psi \; + \; \tfrac{r^{p+1}}{p+1}\, d \psi \ .
\end{equation}

The next lemma is our main technical tool for the classification
of special Killing forms. It states that special Killing forms are 
exactly those forms which translate into parallel forms on the 
metric cone.

\begin{Lem}\label{cone}
Let $(M,\,g)$ be a Riemannian manifold and let $\,\psi\,$ be a $\,p$--form
on $\,M$. Then the associated $\,(p+1)$--form $\,\widehat \psi\,$ on 
the metric cone
$\,\widehat M\,$ is parallel with respect to $\,\widehat \nabla\,$ if and only if
$$
\nabla_X \,\psi 
\;=\; 
\tfrac{1}{p+1}\, X \es d\psi
\qquad\mbox{and}\qquad
\nabla_X \,(d\psi)
\;=\;
-\,(p+1)\, X^* \we \psi
$$
i.e. $\,\widehat \psi\,$ is parallel if and only if $\;\psi$ is a special Killing form 
with constant $\,c=-(p+1)$.
\end{Lem}
\proof
We will first show that a $(p+1)$--form $\hat \psi$ defined on the metric
cone as in (\ref{coneform}) is always parallel in radial direction.
Indeed we have
\bea
{\widehat \nabla}_{\partial_r}\,\psi
&=&
p \,r^{p-1}\,dr \wedge \psi 
\;\;+\;\;
r^p\,dr \wedge {\widehat \nabla}_{\partial_r}\,\psi
\;\;+\;\;
r^p\,d\omega
\;\;+\;\;
\tfrac{r^{p+1}}{p+1}\,{\widehat \nabla}_{\partial_r}\,(d\psi)\\[1.5ex]
&=&
(p \,r^{p-1} - r^p\,\tfrac{p}{r})\,dr \we \psi 
\;\;+\;\;
(r^p  - \tfrac{r^{p+1}}{p+1}\, \tfrac{1}{r}\,(p+1))\,d\psi\\[1.5ex]
&=&
0 \ .
\eea
Next, we compute the covariant derivative of $\,\widehat \psi\,$ in
direction of a horizontal vector field $X$. This yields
\bea
{\widehat \nabla}_{X}\,\psi
&=&
r^p\,{\widehat \nabla}_{X}\,(dr) \wedge \psi
\;\;+\;\;
r^p\,dr \wedge {\widehat \nabla}_{X}\,\psi
\;\;+\;\;
\tfrac{r^{p+1}}{p+1}\, {\widehat \nabla}_{X}\,(d\psi)\\[1.5ex]
&=&
r^{p+1}\,X^* \wedge \psi
\;\;+\;\;
r^p\,dr \wedge \nabla_X \, \psi
\;\;+\;\;
\tfrac{r^{p+1}}{p+1}\,\nabla_X\,(d\psi)
\;\;-\;\;
\tfrac{r^{p}}{p+1}\,dr \wedge (X \es d\psi)\\[1.5ex]
&=&
r^{p+1}\,
\left(
X^* \wedge \psi \;+\; \tfrac{1}{p+1}\,\nabla_X\,(d\psi)
\right)
\;\;+\;\;
r^p\,dr \wedge
\left(
\nabla_X \, \psi \;-\; \tfrac{1}{p+1}\,X \es d\psi
\right) \ .
\eea
From this equation it becomes clear that $\,\widehat \psi\,$ is parallel,
if and only if the two brackets vanish, i.e. if and only if the form 
$\,\psi\,$ on 
$\,M\,$ is
a special Killing form.
\qed

We already know that on Sasakian manifolds, the Killing 1-form  $\,\xi^*\,$ 
together with all forms $\,\xi^* \we (d\xi^*)^k\,$ are special Killing forms. 
As an immediate corollary of Lemma~\ref{cone} we see that a similar statement
is true for all manifolds admitting special  Killing forms of odd degree.
Note that we have to assume the Killing form $\,\psi\,$ to be of odd degree, 
since otherwise $\,d\psi \wedge d\psi = 0\,$  and we could not obtain a new
Killing form.
\begin{Lem}\label{power}
Let $\,\psi\,$ be a special Killing form of odd degree $\,p$, then all
the forms
$$
\psi_k 
\;\;:=\;\;
\psi \we (d\psi)^k
\qquad\qquad 
k=0,\ldots
$$
are special Killing forms of degree $\,p + k(p+1)$.
\end{Lem}
\proof
Let $\,\widehat \psi\,$ be the parallel form associated with
the special Killing form $\,\psi$. Then  the form $\,\widehat \psi_k\,$
associated to $\,\psi_k\,$ turns out to be 
$\,\tfrac{(p+1)^k}{k+1}\,{\widehat \psi}^{\,k+1}$, 
which is again parallel. Hence, $\,\psi_k\,$ is a special 
Killing form. 
\qed

In the proof of the lemma we have used that the power of the
associated form $\,\widehat \psi\,$ is again parallel and can
be written as associated form for some other special
Killing form. The following lemma will show that
this is a general fact, i.e. we have an simple characterization 
of all parallel forms on the metric cone. It turns out that there 
are no other parallel forms on the cone as the ones corresponding 
to special Killing forms  on the base manifold.

\begin{Lem}\label{parconeforms}
Let $\,\omega\,$ be a form on the metric cone $\,{\widehat M}$. Then 
$\,\omega\,$ is parallel with respect to $\,\widehat \nabla\,$ 
if and only if there exists 
a special Killing form
$\,\psi\,$ on $\,M\,$ such that $\,\omega = {\hat \psi}$.
\end{Lem}
\proof
We know already that $\,\widehat \psi\,$ is parallel on the metric cone,
provided that $\,\psi\,$ is a special Killing form on $\,M$. It
remains 
to verify the
opposite direction. Assuming $\,\omega\,$ to be a parallel form on the cone
we write it  as
$$
\omega \;\;=\;\; \omega_0 \;+\; dr \we \omega_1 \ ,
$$
where we consider $\,\omega_0\,$ and $\,\omega_1\,$ as a
$\,r$-dependent 
family
of forms on $\,M$. It is clear that $\,\omega\,$ is parallel in the
radial direction $\,\partial_r\,$ if and only if the same is true for the
two forms $\,\omega_0\,$ and $\,\omega_1$. Let $\,\eta =\eta(r)\,$ be any
horizontal $\,p$-form on $\,\widehat M\,$ considered as family of forms
on $\,M$. Locally we can write
$\,\eta = \sum r^p f_I(r, x)\,dx_{i_1}\wedge \ldots \wedge dx_{i_p}$,
with multi index $\,I=(i_1,\ldots,i_p)$. Then $\,\eta\,$ is parallel
in radial direction if and only if
\bea
0 &=&
\partial_r(r^p f_I(r, x))
\;+\;
r^p f_I(r, x)\,(- \tfrac{p}{r}\,)\\[1.5ex]
&=&
p\,r^{p-1}f_I(r, x)
\;+\;
r^p\,\partial_r(f_I(r, x))
\;-\;
r^{p-1}p\,f_I(r, x)\\[1.5ex]
&=&
r^p\,\partial_r(f_I(r, x)) \ . 
\eea
It follows that $\,f_I(r, x)\,$ does not depend on $\,r$. Hence, we can 
write $\,\eta= r^p\,\eta_0$, where $\,\eta_0$ is a $\,p$-form on $\,M$. 
In particular, we have $\,\omega_0 =r^{p+1}\omega_0^M\,$ and
$\,\omega_1 =r^{p}\,\omega_1^M$, where $\,\omega^M_0\,$ and $\,\omega^M_1\,$
are forms on $\,M$. Next, we consider the covariant derivative of 
the parallel form $\,\omega\,$ in direction of a horizontal vector 
field $\,X$. Here we obtain
\bea
{\widehat \nabla}_X\,\omega
&=&
r^{p+1}{\widehat \nabla}_X\,\omega_0^M
\;\;+\;\;
r^{p+1}\,X^* \we \omega^M_1
\;\;+\;\;
r^pdr \we {\widehat \nabla}_X\,\omega^M_1\\[1.5ex]
&=&
r^{p+1}
\left(
\nabla_X\,\omega_0^M
\;-\;
\tfrac{1}{r}\,dr \we (X \es \,\omega_0^M)
\right)
\\[1.2ex]
&&
\qquad\qquad
\;\;+\;
r^{p+1}X^* \we \omega^M_1
\;\;+\;\;
r^p dr \we \nabla_X\,\omega^M_1 \ .
\eea
From this we conclude that the form 
$\,\omega = r^p\,dr \wedge \omega^M_1 + r^{p+1}\omega_0^M \,$ is
parallel if and only if the following two equations are satisfied
for all vector fields $\,X\,$ on $\,M$
\begin{equation}\label{conepar}
\nabla_X\,\omega^M_1 
\;=\;
X \es \omega^M_0
\qquad\quad\mbox{and}\qquad\quad
\nabla_X\,\omega^M_0 
\;=\;
-\,X^* \we \omega^M_1 \ .
\end{equation}
Using these equations we immediately find:
\bea
d \,\omega^M_0 &=& 0 \;=\; d^*\omega^M_1,\quad
d\,\omega^M_1 = (p+1)\,\omega^M_0,\quad 
d^*\omega^M_0=(n-p)\omega^M_0 \ .\quad 
\eea
In particular, we have $\,\Delta\,\omega^M_1=(p+1)(n-p)\,\omega^M_1\,$
and it is clear that $\,\omega = {\widehat \psi}\,$ for the special
Killing $\,p$-form $\,\psi=\omega^M_1$.
\qed

We have seen that the map $\,\psi \mapsto {\widehat \psi}\,$
defines a 1-1-correspondence between special Killing $\,p$-forms
on $\,M\,$ and parallel $\,(p+1)$-forms on the metric cone $\,\widehat M$.
We will use this fact to describe manifolds admitting
special Killing forms. Let $\,M\,$ be a compact oriented simply connected 
manifold, then the metric cone $\,\widehat M\,$ is either
flat, and the manifold $\,M\,$ has to be isometric to the standard 
sphere, or the cone is irreducible (c.f.~\cite{baer} or \cite{galot1}).
In the latter case we know from the holonomy theorem of M.~Berger that
$\,\widehat M\,$ is either symmetric or its holonomy is one of the
the following groups: 
$\,\SO(m),\,\Sp(m)\cdot \Sp(1),\,\U(m),\,\SU(m),\,\Sp(m),\,
\G_2\,$ or $\,\Spin_7$. An irreducible symmetric space as well as a manifold 
with holonomy $\,\Sp(m)\cdot \Sp(1)\,$ is automatically Einstein 
(c.f.~\cite{besse}).
But it follows from the O'Neill formulas applied to the 
cone, that $\,{\widehat \Ric}(\partial_r,\,\partial_r)=0$,~i.e. the
metric 
cone
can only be Einstein if it is
Ricci-flat. In this case the symmetric space has to be flat and the 
holonomy $\,\Sp(m)\cdot \Sp(1)\,$ restricts further to $\,\Sp(m)\,$ 
(this again can be found in c.f.~\cite{besse}). 

Let $(M,\,g)$ be a compact oriented simply connected manifold not isometric
to the sphere. If $\,\psi\,$ is a special Killing form on $\,M\,$ then
the 
metric
cone $\,\widehat M\,$ is an irreducible manifold with a parallel form
$\,\widehat \psi$. Since any parallel form induces a holonomy reduction,
we see that the above list of possible holonomies is further reduced to
$\U(m),\,\SU(m),\,\Sp(m),\,\G_2,\,$ or $\,\Spin_7$. We will now go
through 
this
list and determine what are the possible parallel forms and how they 
translate
into special Killing forms on $\,M$. The description of possible parallel
forms can be found in \cite{besse}, with the only exception of holonomy
$\,\Sp(m)$. Nevertheless, in this case the parallel forms can be described using the
realization of $\,\Sp(m)$--representation due to H.~Weyl (the result is
also contained in \cite{fukami}). Concerning the translation from special
holonomy on $\,\widehat M\,$ to special geometric structures on $\,M$ we
refer to \cite{baer}, where the explicit constructions are described.

The first case,~i.e. holonomy $\,\U(m)$, is equivalent to $\,\widehat
M\,$ 
being
a K{\"a}hler manifold. In this case all parallel forms are linear combinations
of powers of the K{\"a}hler form. On the other hand, it is well-known that
$\widehat M$ is K{\"a}hler, if and only if $M$ is a Sasakian manifold. 
If the Killing
vector field $\xi$ defines the Sasakian structure on $M$, then 
${\widehat \xi} = r dr \wedge \xi^* + \frac{r^2}{2}d\xi^*$ defines the
K{\"a}hler form on $\widehat M$. Hence, all special Killing forms on a
Sasakian manifold are spanned by the forms $\,\omega_k\,$ given in 
Proposition~\ref{powers}, and they all correspond to the powers of the
K{\"a}hler form on $\widehat M$.

In the next case, $\widehat M$ has holonomy $\SU(m)$ and equivalently is
Ricci-flat and K{\"a}hler. In this situation, there are two additional parallel
forms given by the complex volume form and its conjugate. As real forms
we obtain the real part resp. the imaginary part of the complex volume form.
Because of the O'Neill formulas, the cone is Ricci-flat, if and only
if 
the base manifold 
is Einstein,~i.e. in this case our manifold is Einstein-Sasakian. As special
Killing forms we have the forms $\,\omega_k\,$ and two additional 
forms of degree $\,m$, which can also be described using the 
Killing spinors of an Einstein-Sasakian manifold.

In the third case, $\,\widehat M\,$ has holonomy $\,\Sp(m)\,$ and is 
by definition 
a hyper-K{\"a}hler manifold,~i.e. there are three K{\"a}hler forms
compatible with the metric and such that the corresponding complex 
structures satisfy the quaternionic relations. Here, all
parallel forms are linear combinations of wedge products of powers
of the three K{\"a}hler forms (c.f.~\cite{fukami}). The metric cone is hyper-K{\"a}hler if and only 
if the base 
manifold has a 3-Sasakian structure and the possible special Killing
forms are described by
\begin{Pro}\label{3sasaki}
Let $(M,\,g)$ be a manifold with a 3-Sasakian structure
defined by the Killing 1--forms $\,\eta_1,\,\eta_2\,$ and $\,\eta_3$. Then
all special Killing forms on $\,M\,$ are linear combinations of the
forms $\,\psi_{a,\, b,\, c}\,$ defined for any integers $\,(a, b, c)$  by
\bea
\firstline{.5cm}{\psi_{a,\, b, \,c}
:=
\tfrac{a}{a+b+c}\,
[\eta_1 \we (d\eta_1)^{a-1}] \we (d\eta_2)^b \we (d\eta_3)^c }
\\[1.5ex]
&&\qquad\qquad
+\;\;
\tfrac{b}{a+b+c}\,
(d\eta_1)^a \we [\eta_2 \we (d\eta_2)^{b-1}] \we (d\eta_3)^c 
\\[1.5ex]
&&\qquad\qquad\qquad\qquad
+\;\;
\tfrac{c}{a+b+c}\,
(d\eta_1)^a \we (d\eta_2)^b \we [\eta_3 \we (d\eta_3)^{c-1}] \ .
\eea
\end{Pro}
\proof
Let $\,\phi_i\,$ be the parallel 2-form associated with the Sasakian
structure $\,\eta_i$, for $\,i=1,2,3$, i.e.
$$
\phi_i \;=\;r\,dr \we \eta_i \;\;+\;\;\tfrac{r^2}{2}\,d\eta_i \ .
$$
Then it follows from a simple computation that 
$\,\phi^a_1 \we \phi^b_2 \we \phi^c_3\,$ is a parallel form which
is, up to a factor, associated to the form $\,\psi_{a, b, c}\,$
defined above.
\qed

Next, we have to consider the two exceptional holonomies 
$\,\G_2\,$ resp. $\,\Spin_7$. These holonomies are defined by the existence
of a parallel 3-- resp. 4--form $\psi$ and the only non-trivial parallel 
forms on such a manifold are the linear combinations of $\,\psi\,$ and 
$\,\ast \psi$.
The metric cone has holonomy $\,\G_2$ if and only if the base manifold is a
6-dimensional nearly K{\"a}hler manifold. Here, the parallel 3-form $\,\psi\,$
translates into the K{\"a}hler form $\,\omega\,$ and the parallel 
4-form $\,\ast \psi\,$ translates, up to a constant, into the 3-form 
$\,\ast d \omega$.
To make this more precise, we note the following simple fact
\begin{Lem}
Let $\,\omega\,$ be a $\,p$-form on $\,M\,$ considered as $\,p$-form on the
metric cone $\,\widehat M$. Then the Hodge star operators of $\,M\,$ and
$\,\widehat M\,$ are related by
$$
\ast_{ \widehat M}\,\omega
\;\;=\;\;
r^{n-2p} (\ast_M \omega) \we dr \ .
$$
\end{Lem}
Now, back to the nearly K{\"a}hler case, let $\,\psi = r^2 dr \wedge \omega + 
\tfrac{r^3}{3}d\omega\,$ be the parallel 3-form associated with the K{\"a}hler 
form $\,\omega$.  As in the proof of Lemma~\ref{parconeforms} we conclude
$\,\Delta \omega = 12 \omega$.
Hence, the scalar curvature $\,s_M\,$ of the 6-dimensional nearly 
K{\"a}hler manifold
is normalized to  $\,s_M = 30$. Applying the lemma above yields
\bea
\ast_{\widehat M}\,\psi
&=&
r^2\,\ast_{\widehat M} (dr \we \omega)
\;+\;
\tfrac{r^3}{3}\,\ast_{\widehat M} (d\omega)
\\[1.5ex]
&=&
r^2\,\partial_r \es (\ast_{\widehat M}\,\omega)
\;+\;
\tfrac{r^3}{3}\,\ast_{\widehat M} (d\omega)
\;\;=\;\;
r^4\ast_M\omega 
\;+\;
\tfrac{r^3}{3}\,(\ast_{M} \,d\omega) \we dr \ .
\eea
Since $\,\Delta\, \omega = 12 \,\omega\,$ and $\,d^*\omega = 0$ it follows
$\,d^*d\omega = -\ast_Md\ast_Md\omega =12 \omega\,$ and we obtain
$\,d(\ast_Md\omega) = -12\,\ast_M\,\omega$. Substituting this into the
equation for $\,\ast_{\widehat M}\,\psi$, we find
$$
\ast_{\widehat M}\,\psi
\;\;=\;\;
-\,\tfrac{r^4}{12}\,d(\ast_Md\omega)
\;-\;
\tfrac{r^3}{3}\,
dr \we (\ast_M d\omega)\ .
$$
From where we conclude that $\,\ast_M d\omega\,$ is the special
Killing form on the nearly K{\"a}hler manifold $\,M\,$ corresponding to
the parallel 4-form $\,-3\ast_{\widehat M} \psi\,$ on $\,\widehat M$.

Finally we have to consider the case of holonomy $\Spin_7$. The metric 
cone has holonomy
$\Spin_7$ if and only if $M$ is a 7-dimensional manifold with a weak 
$G_2$-structure.
Here the parallel 4-form $\psi$ on the cone is self-dual, i.e. 
$\ast \psi = \psi$,
and the corresponding special Killing form is just the 3-form defining
the weak $G_2$-structure.

Summarizing our description of compact manifolds with special Killing
forms we have the following
\begin{The}\label{classe}
Let $(M^n,\,g)$ be a compact, simply connected manifold admitting
a special Killing form. Then $\,M\,$ is either isometric to $\,S^n\,$ or
$\,M\,$ is a Sasakian, 3-Sasakian, nearly K{\"a}hler or weak $\,G_2$--manifold.
Moreover, on these manifolds any special Killing form is a linear
combination of the Killing forms described above.
\end{The}

\bigskip


\section{The dimension bound}\label{dimension}

It is well-known and easy to check that for twistor operator $\,T\,$ the
operator $\,T^*T\,$ is elliptic. 
Hence, the space of conformal Killing forms is finite dimensional on compact manifolds. 
However, in this section we will prove that the space of conformal Killing forms is 
finite dimensional on any connected manifold. More precisely,  we have

\begin{The}\label{dimbound}
Let $\,(M,\,g)\,$ be an $\,n$-dimensional connected Riemannian manifold and 
denote with $\,{\mathcal CK}^p(M)\,$ the space of conformal Killing p-forms, then
$$
\dim \,{\mathcal CK}^p(M)
\;\;\le\;\;
{n+2 \choose p+1}
$$
with equality attained on the standard sphere. Moreover, if a manifold 
admits the maximal possible number of linear independent conformal
Killing $p$-forms, with $\,1<p<n-1$, then it is conformally flat.
\end{The}

The idea of the proof is to  construct a vector bundle together with a connection, called {\it Killing connection}, such that conformal Killing forms are in a 1-1-correspondence to parallel sections for this connection. It then follows immediately that the dimension of the space of conformal Killing forms is bounded 
by the rank of the constructed vector bundle. 

By definition, the covariant derivative of a conformal Killing $p$--form $\,\psi\,$ involves $\,d\psi\,$ and $\,d^*\psi$. Computing the covariant derivatives of $\,d\psi\,$ and $\,d^*\psi\,$  we obtain an expression involving only $\,\psi\,$ and $\,d d^*\psi$. Finally we have to compute the covariant derivative
of $\,d d^*\psi\,$ which leads to an expression involving only $\,\psi, d\psi\,$ and $\,d^*\psi$. Collecting the covariant derivatives
we can formulate the result of the computations as follows.
Let 
$\,{\hat \psi} := (\psi,\, d\psi,\, d^*\psi,\, d d^*\psi)$, then
$\,{\hat \psi}$ is a section of 
$\,{\mathcal E}^p(M):=\Lambda^p T^*M\oplus\Lambda^{p+1} T^*M\oplus\Lambda^{p-1} 
T^*M\oplus \Lambda^p T^*M\,$ and we have 
$\,\nabla_X {\hat \psi} = A(X) \,{\hat \psi}$, where $\,A(X)\,$ is a certain $4\times 4$-matrix
with coefficients which are endomorphisms of the form bundle depending
on the vector field $\,X$. Here the components of $\,\nabla_X {\hat \psi}\,$ are
the covariant derivatives of the components of $\,{\hat \psi}$.
The Killing connection $\widetilde \nabla$ is then a connection on $\,{\mathcal E}^p(M)\,$
defined as $\,{\widetilde \nabla}_X := \nabla_X - A(X)\,$ and the conformal Killing forms
are by definition the first components of parallel sections of $\,{\mathcal E}^p(M)$.
Hence, the rank of the bundle ${\mathcal E}^p(M)$ is an upper bound on the dimension of the space of conformal Killing forms,~i.e.
$$
\dim\,{\mathcal CK}^p(M)
\; \le \;
2\,{n \choose p}
\;+\;
{n \choose p-1 }
\;+\;
{n \choose p+1}
\;=\;
{n + 1 \choose p}
\;+\;
{n + 1 \choose p+1}
\;=\;
{n+2 \choose p+1} \ .
$$

It follows from Proposition~\ref{sphere} that this upper bound is attained on
the standard sphere. Moreover, if on $M$ exists the maximal possible number of linearly independent conformal Killing forms then the  map  $\,{\mathcal E}^p(M)\rightarrow \Lambda^p(T^*_xM)$, defined as projection onto the first component and evaluation
in the point $x$, is obviously surjective,~i.e. any p-form in $\,\Lambda^p(T^*_xM)\,$
can be extended to a conformal Killing form. In this situation, and with 
$\,1<p<n-1$,  a curvature
calculation shows that the manifold has to be conformally flat (c.f.~\cite{ka}).

In the remaining part of this section we will show the existence of the 
Killing connection, which then concludes the proof the Theorem~\ref{dimbound}.
The covariant derivatives of $\,d\psi\,$ and $\,d^*\psi\,$ for a conformal 
Killing form $\,\psi\,$ can be obtained by a direct calculation starting
from the definition. In order to give the explicit formulas we introduce
the notation $\;R^+(X)\psi:=\sum\,e_j\,\wedge\,R_{X,e_j}\psi\,$ and 
$\,R^-(X)\psi:=\sum\,e_j\es R_{X,\,e_j}\psi$, where $\,\{e_i\}\,$ is a 
local ortho-normal basis  and  $\psi$ is any differential form.

\begin{Pro}\label{nabla-d}
Let $\psi$ be a conformal Killing p-form, then for all vector fields $X$
\bea
\nabla_X(d\,\psi) &=& \quad
\tfrac{p+1}{p}\,R^+(X)\,\psi\;+\;
\tfrac{p+1}{p(n-p+1)}\,X\,\wedge\,d\,d^*\,\psi \\[1.5ex]
\nabla_X(d^*\,\psi)
&=&
-\,\tfrac{n-p+1}{n-p}\,R^-(X)\,\psi\;+\;
\tfrac{1}{p}\,X\es d\,d^*\,\psi\;-\; 
\tfrac{n-p+1}{p(n-p)}\,X\,\lrcorner\,q(R)\,\psi \ .
\eea
\end{Pro}

It remains to show that the covariant derivative of $\,d d^*\psi\,$  is
an expression only involving  $\,\psi, d\psi\,$ and
$\,d^*\psi$. With regard to equation~(\ref{covariant}) it suffices to 
consider $\,T(dd^*\psi)$. For this we derive  Weitzenb\"ock formulas which 
then also imply Proposition~\ref{nabla-d}. The twistor operator $T$ was 
defined as the composition of the covariant derivative $\nabla$ with the 
projection $\pr_{\Lambda^{p,1}}$. Similarly we obtain operators like
$\,Td\psi\,$ or $\,Td^*\psi$ by applying certain projections to 
$\,\nabla^2 \psi$. Hence, it suffices to consider relations between such projections
which then translate into Weitzenb\"ock formulas for the corresponding
differential operators. As a first projection we define
\bea
\pr^+_1 \;:\; 
T^*M \otimes T^*M \otimes \Lambda^pT^*M 
& \rightarrow & 
T^*M \otimes \Lambda^{p+1}T^*M
\quad  \rightarrow  \quad 
\Lambda^{p+1,1}T^*M 
\\[1.5ex]
e_1 \otimes e_2 \otimes \psi 
&\mapsto &
e_1 \otimes (e_2 \we \psi)
\quad \mapsto \quad
\pr_{\Lambda^{p+1,1}}\,(e_1 \otimes (e_2 \we \psi)) \ .
\eea
Let $\,\psi\,$ be any p-form then $\,\nabla^2\psi\,$ is a section of 
$\,T^*M \otimes T^*M \otimes \Lambda^pT^*M\,$ and it is easy to show that
$\, \pr^+_1\,(\nabla^2\psi) = T(d\psi)$. Next we need the map
\bea
\pr^+_2 \;:\; 
T^*M \otimes T^*M \otimes \Lambda^pT^*M 
& \rightarrow & 
T^*M \otimes \Lambda^{p,1}T^*M
\quad  \rightarrow  \quad 
\Lambda^{p+1,1}T^*M 
\\[1.5ex]
e_1 \otimes e_2 \otimes \psi 
& \mapsto & 
e_1 \otimes \pr_{\Lambda^{p,1}}(e_2 \otimes \psi)
\mapsto 
\pr_{\Lambda^{p+1,1}}(e_1  \we \pr_{\Lambda^{p,1}}(e_2 \otimes \psi)  ) .
\eea
In this case there appears a new first order differential operator, which
we denote by $\,\theta^+$. It maps sections of $\, \Lambda^{p,1}T^*M\,$
into sections of $\,\Lambda^{p+1,1}T^*M\,$
and is defined as $\,\pr \circ \nabla$, where $\,\pr\,$ is the projection 
$\,T^*M \otimes \Lambda^{p,1}T^*M \rightarrow  \Lambda^{p+1,1}T^*M\,$
defined above (as the second map in the definition of $\pr^+_2$). We have
$$
\, \pr^+_2\,(\nabla^2\psi)\;\; =\;\; \theta^+ T(\psi) \ .
$$
Then we need a third projection which will produce the
curvature term. We define it as
\bea
&&
\pi^+ \;:\; 
T^*M \otimes T^*M \otimes \Lambda^pT^*M 
\quad\rightarrow \quad
\Lambda^2 T^*M \otimes \Lambda^{p}T^*M
\quad \rightarrow  \quad 
\Lambda^{p+1,1}T^*M 
\\[1.5ex]
&&\qquad
e_1 \otimes e_2 \otimes \psi 
\;\mapsto \;
(e_1 \we e_2) \otimes \psi
\;\mapsto \;
\pr_{\Lambda^{p+1,1}}\,
(\sum\,  e_i \es (e_1 \we e_2) \otimes (e_1 \we \psi) ) \ .
\eea
Let $\psi$ be any p-form, then $\,\nabla^2\psi\,$ is a section of
$\,T^*M \otimes T^*M \otimes \Lambda^pT^*M\,$ and the first map
in the definition of $\,\pi^+\,$ maps this section to the
curvature $\,R(\cdot, \cdot)\psi$. Computing the result of the
second map we obtain
$$
\pi^+\,(\nabla^2\psi)
\;=\;
-\,\tfrac{1}{p}\,R^+(\cdot)\,\psi
\;-\;
\tfrac{1}{p(n-p)}\; . \we q(R)\,\psi \ .
$$
Having defined these three projections it is an elementary calculation to 
prove that they satisfy the following linear relation
$$
(p+1)\,\pi^+
\;\;+\;\;\;
\pr^+_1
\;\;\;=\;\;\;
\tfrac{p+1}{p}\,\pr^+_2 \ .
$$
To obtain a twistor Weitzenb\"ock formula we only have to apply this
relation to $\,\nabla^2\psi\,$ and to substitute the expressions for the
three different projections of $\nabla^2\psi$. The result is
\begin{Lem}\label{twistor1}
Let $\,\psi\,$ be any p-form then:
\bea
T(d \psi)
&=&
\tfrac{p+1}{p}\,\theta^+( T \psi)
\;+\;
\tfrac{p+1}{p}
\,R^+(\cdot)\,\psi
\;+\;
\tfrac{p+1}{p(n-p)}\, \cdot \we q(R)\,\psi \ .
\eea
\end{Lem}

By defining similar projections or by applying the Hodge star operator to the equation
of Lemma~\ref{twistor1}, with $\psi$ replaced by $\ast \psi$, we obtain a 
corresponding formula for $\,T(d^* \psi)$. Here appears an operator $\theta^-$, 
which is defined as $\theta^+$, only with the wedge product replaced by the contraction. 
In this case the result is

\begin{Lem}\label{twistor2}
Let $\psi$ be any p-form then:
\bea
T(d^* \psi)
&=&
\tfrac{n-p+1}{n-p}\,\theta^-( T \psi)
\;-\;
\tfrac{n-p+1}{n-p}
\,R^-(\cdot)\,\psi
\;-\;
\tfrac{n-p+1}{p(n-p)}\, \cdot \es q(R)\,\psi \ .
\eea
\end{Lem}

If $\psi$ is a conformal Killing form then $T\psi=0$ and the summands
with $\theta^+$ resp. $\theta^-$ vanish. Substituting the expressions
for $Td\psi$ resp. $Td^*\psi$ into equation~(\ref{covariant}) proves
Proposition~\ref{nabla-d}. 

Finally we have  to show that the covariant derivative of $dd^*\psi$
for a conformal Killing p-form $\psi$ can be obtained from 
$d\psi$ resp. $d^*\psi$ by applying certain bundle homomorphisms.
We replace in the formula of Proposition~\ref{twistor1} the $p$-form $\psi$ with
$d^*\psi$  and $p$ with $p-1$ to obtain
$$
T(d d^* \psi)
\;=\;
\tfrac{p}{p-1}\,\theta^+( T d^*\psi)
\;+\;
\tfrac{p}{p-1}
\,R^+(\cdot)\,d^*\psi
\;+\;
\tfrac{p}{(p-1)(n-p+1)}\, \cdot \we q(R)\,d^*\psi \ .
$$
It remains to investigate the summand with $\theta^+( T d^*\psi)$.
Since $\psi$ is a conformal Killing form we can use Lemma~\ref{twistor2}
to replace $Td^*\psi$,~i.e. we have 
\bea
\firstline{.5cm}{
\theta^+( T d^*\psi)
\;=\;
c_1\sum\,\theta^+(e_i \otimes R^-(e_i)\,\psi)
\;+\;
c_2\sum\,\theta^+(e_i \otimes e_i \es q(R)\psi)}\\[1.5ex]
&=&
c_1\sum\,\pr (e_j \otimes e_i \otimes \nabla_{e_j}(R^-(e_i)\,\psi))
\;+\;
c_2\sum\,\pr (e_j \otimes e_i \otimes e_i \es \nabla_{e_j}(q(R)\psi)) \ ,
\eea
where the constants $c_1$ and $c_2$ are given by Lemma~\ref{twistor2} and
where we do the calculation in a point with $\nabla e_i = 0$. Computing the 
covariant derivative of $R^-(\cdot)\psi$ and $q(R)\psi$ easily leads to
\begin{Lem}
Let $\psi$ be any differential form then for any vector fields $X$, $Y$
\bea
(i)\qquad
\nabla_X(R^-(Y)\psi)
&=&
R^-(\nabla_XY)\psi
\;+\;
R^-(Y)\nabla_X\psi
\;+\;
(\nabla_X R)^-(Y)\psi
\\[1.5ex]
(ii)\qquad
\nabla_X(q(R)\psi)
&=&
q(\nabla_XR)\psi \;+\; q(R)\nabla_X\psi
\eea
\end{Lem}

Using this lemma we can substitute the summands $\nabla_{e_j}(R^-(e_i)\,\psi)$
and $\nabla_{e_j}(q(R)\psi))$ in the formula for $\theta^+( T d^*\psi)$
and see that it indeed only involves the covariant derivative of $\psi$, which
for the conformal Killing form $\psi$ is an expression in $d\psi$ and $d^*\psi$.
Summarizing the calculations we see that the covariant derivative of the
four sections $\psi, d\psi, d^*\psi$ and $dd^*\psi$ can be obtained from
these sections by applying certain bundle homomorphisms. Hence we can collect
the covariant derivatives to define a Killing connection (as explained above),
which then concludes the proof of the dimension bound.


\section{Further results}

In this section we state (without proof) several further results on
conformal Killing forms. We start with compact K\"ahler
manifolds, where it is easy to show that any Killing form has to be parallel
(c.f.~\cite{yama3}). More generally we proved in~\cite{au} 

\begin{The}
On a compact
K{\"a}hler manifold $M^{2m}$ any conformal Killing form $\psi$ has to be of the
form
$$
\psi
\;=\;
L^{k-1}\phi
\;+\;
L^kf
\;+\;
\psi_0 \ ,
$$
where $L$ denotes the wedging with the K\"ahler form, $\phi$ is a special
2-form with associated function $f$ and $\psi_0$ is any parallel form.
Conversely any special 2-forms defines in this way conformal Killing forms
on $M$ in any even degree.
\end{The}

Special 2-forms are defined as primitive (1,1)-forms satisfying an
additional differential equation. In particular, it follows for a
special 2-form $\phi\,$ that $Jd^*\phi$ is exact, thus defining the function $f$
(up to constants). Special 2-forms are closely related to Hamiltonian
2-forms, which were studied and locally classified in~\cite{acg}. In particular,
if $m>2$ then any special 2-form is the primitive part of a Hamiltonian
2-form and vice versa. Starting from the differential of eigenfunctions
for the minimal eigenvalue of the Laplacian on the complex projective space
one easily can construct special 2-forms. Hence, the complex projective
space admits conformal Killing forms in any even degree. Besides the
complex projective spaces there are several other examples of compact
K\"ahler manifolds with conformal Killing forms. This is in contrast to
results in~\cite{yama}. However, it turns out that the proofs in~\cite{yama} 
contain serious gaps.

\medskip

Let $(M^n, g)$ be a Riemannian manifold such that the holonomy group of
$M$ is a proper subgroup of $O(n)$. In this situation the bundle of
forms decomposes into a sum of parallel subbundles, which are preserved
by the Laplace operator $\Delta$ and the curvature endomorphism $q(R)$.
For any form $\psi$ we have the corresponding holonomy decomposition 
$\psi= \sum \psi_i$, where the forms $\psi_i$ are the projections of $\psi$ onto the 
parallel subbundles. We would like to use the characterization of Killing 
forms given in Corollary~\ref{integrab} to conclude that a form $\psi$, with holonomy decomposition $\psi= \sum \psi_i$, is a Killing form if and only if all components
$\psi_i$ are Killing forms. This not true in general since the components
$\psi_i$ of a coclosed form $\psi$ need not to be coclosed. However, the
statement is true for Killing $m$-forms on a 2m-dimensional manifold and
for Killing forms on manifolds with $G_2$- resp. $Spin_7$-holonomy.
In the case of compact manifolds with holonomy $G_2$ (and similarly for
manifolds with holonomy $Spin_7$, we can derive the following result.

\begin{The}\label{theg2}
Let $(M^7,\,g)$ be a compact manifold with holonomy  $G_2$. Then 
any Killing form and any $\ast$--Killing form is parallel.
Moreover, any conformal Killing $p$--form, with $\,p\ne 3,4$, is parallel.
\end{The}

First of all we note that on a compact Ricci-flat manifold any conformal 
vector field  has to be a Killing vector field and any Killing
vector field has to be parallel. This follows from results of M.~Obata
in~\cite{ob} and Corollary~\ref{integrab}. Hence, on a compact manifold 
with holonomy $G_2$ or $Spin_7$ any conformal Killing 1-form has to be parallel.
Moreover, it is easy to show that on an Einstein manifold any conformal
Killing form $\psi$ is either coclosed,~i.e. Killing, or $d^*\psi$ is
a non-trivial Killing vector field. Thus we obtain that any 
conformal Killing $p$--form, with $\,p\ne 3,4$, is either closed or
coclosed, with a similar statement for $Spin_7$-manifolds.
In the end it remains to consider Killing resp. $\ast$-Killing
forms lying in one of the parallel subbundles of the 2- resp. 3-form
bundle of a $G_2$-manifold. Using additional twistor operators and
explicit formulas for the projections onto the parallel subbundles
it is easy to derive a contradiction to the norm estimate of
Lemma~\ref{kernel}, which proves that any Killing resp. $\ast$-Killing
form has to be parallel. 

\medskip

A special case of a manifold with restricted holonomy is a Riemannian
product $M = M_1 \times M_2$. In this case the holonomy decomposition of the 
form bundle coincides with the decomposition 
$\,\Lambda^p(T^*M) = \sum^p_{r=0} \Lambda^r(T^*M_1)\otimes\Lambda^{p-r}(T^*M_2)\,$
and it is easy to verify that a form $\psi$ is a Killing form if 
and only if all its components $\psi_i$ are Killing forms. More
generally we can show (c.f.~\cite{au2}) the following

\begin{The}
  Every conformal Killing form on a Riemannian product $M=M_1\times M_2$
  is a sum of forms of the following types: parallel forms,
  pull--backs of Killing forms on $M_1$ or $M_2$, and wedge products of
  the volume form of $M_1$ (or $M_2$) with the pull--back of a 
  $\ast$-Killing form on $M_2$ (resp. $M_1$).
\end{The}

Finally we want to give a new description of a curvature condition, which already 
appears  in~\cite{ka}. Using the notation of Section~\ref{dimension} the condition 
can be reformulated as

\begin{Pro}\label{curvature1}
Let $(M^n,\,g)$ be a Riemannian manifold with a conformal Killing
p--form $\,\psi$, then for any vector fields $\,X,\,Y\,$ the following
equation is satisfied:
\bea
\firstline{.5cm}{R(X,\,Y)\,\psi
\;=\;
\tfrac{1}{p(n-p)}
\left(
 Y \wedge X \es \, - \, X \wedge Y \es 
\right) 
q(R)\,\psi }
\\[1.5ex]
&&\quad
\;-\;
\tfrac{1}{p}
\left(
X \es R^+(Y) \, - \, Y \es R^+(X)
\right)
\psi  
\;-\;
\tfrac{1}{n-p}
\left(
X \we R^-(Y) \, - \, Y \we R^-(X)
\right)
\psi \ .
\eea
\end{Pro}

The proof of this proposition is simple local calculation, which is contained
in the computation of the components of the Killing connection.
Considering $\,R(\cdot, \cdot)\,\psi \,$ as a section of 
$\,\Lambda^p(T^*M) \otimes \Lambda^2(T^*M)$, we can write the above curvature 
condition in a much shorter form. Indeed, we have a decomposition of 
the tensor product  $\,\Lambda^p(T^*M) \otimes \Lambda^2(T^*M)\,$ corresponding
to the following isomorphism of 
$O(n)$--representations:
\be \label{deco3}
\Lambda^p V^* \,\otimes\, \Lambda^2 V^*
\;\cong\;
\Lambda^p V^* \oplus\Lambda^{p+1,1} V^* \oplus
\Lambda^{p-1,1} V^* \oplus\Lambda^{p+2} V^*
\oplus\Lambda^{p-2}V^*\oplus\Lambda^{p,2}V^* \ .
\ee
Here $\,\Lambda^{p,2}V^*\,$ is defined as the irreducible representation
which has as highest weight the sum of the highest weights
of $\,\Lambda^p V^*\,$ and  $ \,\Lambda^2 V^*$. 
It is easy to find explicit expressions for the projections onto the six
summands on the right hand side of~(\ref{deco3}), denoted as
$\,\pr_{\Lambda^p},\,\pr_{\Lambda^{p\pm1,1}},\,\pr_{\Lambda^{p\pm 2}}\,$
and $\,\pr_{\Lambda^{p,2}}$. It then follows that the 
projections of $\,R\,(\cdot, \cdot)\,\psi\,$ onto the summands 
$\,\Lambda^{p\pm 2} T^*M\,$ vanish because of the Bianchi identity and that the 
projection of $\,R\,(\cdot, \cdot)\,\psi\,$ onto $\,\Lambda^{p} T^*M\,$ is 
precisely $\,q(R)\psi$. Moreover, it is also not difficult to show that the 
curvature relation can be written as

\begin{Cor}\label{curvature4}
Let $(M^n,\,g)$ be a Riemannian manifold with a conformal Killing
form $\,\psi$. Then
\be 
\pr_{\Lambda^{p,2}}(\,R\,(\cdot, \cdot)\,\psi\,) \;\;=\;\;0 \ .
\ee

If $\,\psi\,$ is coclosed, then the additional equation
$\,\pr_{\Lambda^{p-1,1}}(\,R\,(\cdot, \cdot)\,\psi\,) = 0\,$  is satisfied. Similarly,
if $\,\psi\,$ is closed then the additional equation 
$\,\pr_{\Lambda^{p+1,1}}(\,R\,(\cdot, \cdot)\,\psi\,) = 0\,$
holds.
\end{Cor}

We note that it is possible to give an alternative proof of the
curvature condition of Proposition~\ref{curvature1}, using the Killing connection. 
Indeed, since a conformal Killing form is parallel with respect to the Killing 
connection,  it follows that the curvature  of the Killing connection 
applied to a conformal Killing form has to vanish. This yields four
equations corresponding to the four components of the bundle
$\,{\mathcal E}^p(M)$. The first of these equations turns out to be
equivalent to the curvature condition of Proposition~\ref{curvature1}.


\begin{appendix}
\setcounter{equation}{0}
\section{Proof of Proposition~\ref{twistor}}

In this appendix we will give the proof of Proposition~\ref{twistor}.
Let $\psi$ be a {\it twistor spinor}, i.e. a spinor satisfying for
all vector fields $X$ the equation
$
\nabla_X \psi \;=\; -\,\frac{1}{n}\,X \cdot D\psi ,
$
with Clifford multiplication $\,\cdot\,$ and Dirac operator $D$. Given two such 
twistor spinors, $\psi_1$ and $\psi_2$, we introduced k-forms $\omega_k$, defined
on tangent vectors $X_1, \ldots , X_k$  by
$$
\omega_k(X_1, \ldots, X_k)
:=
\la (X_1 \wedge \ldots \wedge X_k)\cdot \psi_1, \psi_2 \ra \ .
$$
In order to prove that the $\omega_k$'s are indeed twistor forms we
compute first the covariant derivative
$
(\nabla_{X_0}\omega_k)(X_1, \ldots, X_k) \ .
$
Without loss of generality we will do the calculation for a point $p \in M$
and with vector fields $X_i$ satisfying $\nabla_{X_i}X_j=0\,$ in $p$. We obtain
\bea
\firstline{.1cm}
{(\nabla_{X_0}\omega_k)(X_1, \ldots, X_k)
\;=\;
\nabla_{X_0}(\omega_k(X_1, \ldots, X_k))}\\[1.5ex]
&=&
\la  [X_1 \wedge \ldots \wedge X_k]\cdot\nabla_{X_0}( \psi_1), \psi_2 \ra
\;+\;
\la 
[X_1 \wedge \ldots \wedge X_k]\cdot\psi_1,\nabla_{X_0} \psi_2 \ra
\\[1.5ex]
&=&
-\,\tfrac{1}{n}\,
\la  [X_1 \wedge \ldots \wedge X_k]\cdot X_0 \cdot D \psi_1, \psi_2 \ra
\;-\;
\tfrac{1}{n}\,
\la [X_1 \wedge \ldots \wedge X_k]\cdot\psi_1,X_0 \cdot D \psi_2 \ra
\\[1.5ex]
&=&
-\,\tfrac{1}{n}\,
\la  [X_1 \wedge \ldots \wedge X_k]\cdot X_0 \cdot D \psi_1, \psi_2 \ra
\;-\;
\tfrac{1}{n}\,
\epsilon
\,
\la \psi_1, (X_1 \wedge \ldots \wedge X_k)\cdot X_0 \cdot D \psi_2 \ra ,
\eea
where $\epsilon = (-1)^{\tfrac{k(k+1)}{2}}$.
Using the formula 
$\,
\omega \cdot X = (-1)^k ( X \wedge \omega + X \lrcorner \,\omega)
$,
valid for any k-form $\omega$ and any vector field $X$, we can further
reformulate the expression for $\nabla_{X_0}\omega_k$. Setting $X_0 = X_1$
and summing over an orthonormal basis $\{e_i\}$ we find
\bea
\firstline{.1cm}{
d^*\omega_k(X_2, \ldots, X_k)
\;=\;
-\,\sum\,(\nabla{e_i}\,\omega_k)(e_i,X_2, \ldots, X_k )}
\\[1.5ex]
&=&
\tfrac{n-k+1}{n}\,(-1)^{k}\, 
\left(
\la [X_2 \wedge \ldots \wedge X_k]\cdot D\psi_1, \psi_2 \ra \;+\;
\epsilon\,
\la \psi_1, [X_2 \wedge \ldots \wedge X_k] \cdot D \psi_2 \ra 
\right).
\eea
Hence,
\bea
(X_0 \we d^*\omega_k)(X_1, \ldots, X_k)
&=&
\tfrac{n-k+1}{n}\,(-1)^{k}\, 
\la (X_0 \es [X_1 \wedge \ldots \wedge X_k]\cdot D\psi_1, \psi_2 \ra 
\\[1.5ex]
&& 
\qquad \quad
\;+\;
\epsilon\,\tfrac{n-k+1}{n}\,(-1)^{k}\,
\la \psi_1, X_0 \es [X_1 \wedge \ldots \wedge X_k] \cdot D \psi_2 \ra 
.
\eea
A similar calculation for $d\omega_k$  yields
\bea
\firstline{.1cm}{
d\omega_k(X_0, \ldots, X_k)
\;=\;
\sum\,(-1)^i(\nabla_{X_i}\omega_k)(X_0, \dots {\hat X_i} \dots, X_k)}
\\[1.5ex]
&=&
\tfrac{k+1}{n}\,
(-1)^{k+1}\,
(
\la 
[X_0 \wedge \ldots \wedge X_k]
\cdot D\phi_1, \phi_2 
\ra 
\;+\;
\epsilon
\la 
\phi_1, [X_0 \wedge \ldots \wedge X_k] \cdot D\phi_2 
\ra 
)
\eea
Here we used the simple fact that
$\,
\sum\,
(-1)^i\,X_i \es [ X_0 \wedge \ldots \hat{X_i} \ldots \wedge X_k]
= 0 \ .
\,$
Comparing these expressions for $\nabla_{X_0}\omega_k$,
$X_0 \es d\omega_k$ and $X_0 \wedge d^*\omega_k$ we immediately
conclude that $\omega_k$ is a twistor form, i.e. it satisfies
the equation
$$
\nabla_{X_0}\omega_k
\;=\;
\frac{1}{k+1}\,X_0 \es d\omega_k
\;-\;
\frac{1}{n-k+1}\,
X_0 \wedge d^*\omega_k \ .
$$

\end{appendix}


\pagebreak


 \labelsep .5cm


\begin{thebibliography}{22}

 \bibitem{acg}
  \textsc{Apostolov, V., Calderbank, D., Gauduchon, P. },
  \textit{Hamiltonian 2-forms in K{\"a}hler geometry I},
  \textrm{math.DG/0202280 (2002)}.



 \bibitem{baer}
  \textsc{B{\"a}r, C.}
  \textit{Real Killing spinors and holonomy.}
  \textrm{Comm. Math. Phys. 154 (1993), no. 3, 509--521}.


 \bibitem{be1}
  \textsc{Benn, I. M.; Charlton, P.; Kress, J. },
  \textit{Debye potentials for Maxwell and Dirac fields from a generalization 
    of the   Killing-Yano equation}
  \textrm{J. Math. Phys. 38 (1997), no. 9, 4504--4527}.


 \bibitem{be2}
  \textsc{Benn, I. M.; Charlton, P.},
  \textit{Dirac symmetry operators from conformal Killing-Yano tensors}
  \textrm{ Classical Quantum Gravity 14 (1997), no. 5, 1037--1042}.

 \bibitem{besse}
  \textsc{Besse, A.L.}
  \textit{Einstein manifolds}
  \textrm{Ergebnisse der Mathematik und ihrer Grenzgebiete (3);
    Springer-Verlag, Berlin, 1987. }




 \bibitem{branson}
  \textsc{Branson, T.},
  \textit{Stein-Weiss operators and ellipticity}
  \textrm{J. Funct. Anal. 151 (1997), no. 2, 334--383}.



 \bibitem{sala}
  \textsc{Falcitelli, M.; Farinola, A.; Salamon, S. },
  \textit{Almost-Hermitian geometry }
  \textrm{Differential Geom. Appl. 4 (1994), no. 3, 259--282. }



 \bibitem{galicki}
 \textsc{Boyer, C. P.; Galicki, K.; Mann, B. M.},
  \textit{The geometry and topology of $3$-Sasakian manifolds }
  \textrm{ J. Reine Angew. Math. 455 (1994), 183--220.  }



 \bibitem{galot1}
 \textsc{Gallot, S.;},
  \textit{\'{E}quations diff\'{e}rentielles caract\'{e}ristiques de la sph\`{e}re }
  \textrm{ Ann. Sci. \'{E}cole Norm. Sup. (4) 12 (1979), no. 2, 235--267. }






 \bibitem{gray2}
 \textsc{Fernandez, M.;Gray, A.},
  \textit{Riemannian manifolds with structure group $G\sb{2}$}
  \textrm{Ann. Mat. Pura Appl. (4) 132 (1982), 19--45 (1983). }


 \bibitem{fia}
 \textsc{Friedrich, Th.; Kath, I.; Moroianu, A.; Semmelmann, U.},
  \textit{On nearly parallel $G\sb 2$-structures. }
  \textrm{J. Geom. Phys. 23 (1997), no. 3-4, 259--286}.



 \bibitem{fukami}
 \textsc{Fukami, T.},
  \textit{
Invariant tensors under the real representation of symplectic group and their 
applications.  }
  \textrm{Tohoku Math. J. (2) 10 1958 81--90. }




 \bibitem{galot}
 \textsc{Gallot, S.; Meyer, D.  },
  \textit{Op{\'e}rateur de courbure et laplacien des formes 
    diff{\'e}rentielles    d'une vari{\'e}t{\'e} riemannienne. }
  \textrm{ J. Math. Pures Appl. (9) 54 (1975), no. 3, 259--284.}


\bibitem{gray}
 \textsc{Gray, A.},
  \textit{Vector cross products on manifolds. }
  \textrm{Trans. Amer. Math. Soc. 141 (1969), 465--504. }



\bibitem{graywolf}
 \textsc{Gray, A.; Wolf, J.},
  \textit{Homogeneous spaces defined by Lie group automorphisms. I}
  \textrm{ J. Differential Geometry 2 (1968), 77--114 }.








 \bibitem{yama}
  \textsc{Jun, J.-B.; Ayabe, S.; Yamaguchi, S},
  \textit{On the conformal Killing $p$-form in compact Kaehlerian manifolds.}
  \textrm{ Tensor (N.S.) 42 (1985), no. 3, 258--271}.


 \bibitem{ka}
  \textsc{Kashiwada, T.},
  \textit{On conformal Killing tensor.}
  \textrm{Natur. Sci. Rep. Ochanomizu Univ. 19 1968 67--74.}


 \bibitem{ta4}
  \textsc{ Kashiwada, T.; Tachibana, S.},
  \textit{On the integrability of Killing-Yano's equation}
  \textrm{J. Math. Soc. Japan 21 1969 259--265.}



 \bibitem{au}
  \textsc{Moroianu, A.; Semmelmann, U.},
  \textit{Twistor forms on K{\"a}hler manifolds}
  \textrm{ preprint (2002), math.DG/0204322}.


 \bibitem{au2}
  \textsc{Moroianu, A.; Semmelmann, U.},
  \textit{Twistor forms on Riemannian products}
  \textrm{ preprint (2002)}.



 \bibitem{ob}
  \textsc{Obata, M.},
  \textit{The conjectures on conformal transformations of Riemannian manifolds. }
  \textrm{ J. Differential Geometry 6 (1971/72), 247--258.}





 \bibitem{penrose}
  \textsc{Penrose, R.; Walker, M},
  \textit{On quadratic first integrals of the geodesic equations for 
  type $\{22\}$ spacetimes.  }
  \textrm{Comm. Math. Phys. 18 1970 265--274.  }

 \bibitem{ta1}
  \textsc{Tachibana, S.},
  \textit{On Killing tensors in Riemannian manifolds of positive curvature operator}
  \textrm{Tohoku Math. J. (2) 28 (1976), no.2, 177--184.}




 \bibitem{ta2}
  \textsc{Tachibana, S.; Yu, W.N.},
  \textit{On a Riemannian space admitting more than one Sasakian structures.}
  \textrm{Tohoku Math. J. (2) 22 (1970), 536--540}.


 \bibitem{ta3}
  \textsc{Tachibana, S.},
  \textit{On conformal Killing tensor in a Riemannian space}
  \textrm{Tohoku Math. J. (2) 21 1969 56--64.}


 \bibitem{ta5}
  \textsc{Tachibana, S.},
  \textit{On Killing tensors in a Riemannian space.}
  \textrm{Tohoku Math. J. (2) 20 1968 257--264. }


 \bibitem{yama3}
  \textsc{Yamaguchi, S.},
  \textit{On a Killing $p$-form in a compact K{\"a}hlerian manifold}
  \textrm{Tensor (N.S.) 29 (1975), no. 3, 274--276. }

\bibitem{yama3.5}
  \textsc{Yamaguchi, S.},
  \textit{On a horizontal conformal Killing tensor of degree $p$ in a  Sasakian space. }
  \textrm{Ann. Mat. Pura Appl. (4) 94 (1972), 217--230.}


 \bibitem{yama4}
  \textsc{Yamaguchi, S.},
  \textit{On a conformal Killing $p$-form in a compact Sasakian space. }
  \textrm{Ann. Mat. Pura Appl. (4) 94 (1972), 231--245.}




 \bibitem{yano}
  \textsc{Yano, K.},
  \textit{Some remarks on tensor fields and curvature. }
  \textrm{Ann. of Math. (2) 55, (1952). 328--347.}




\end{thebibliography}
\end{document}